%% file: Papier31Aout2023.tex
\documentclass{article}
\usepackage[title,titletoc,toc]{appendix}
\usepackage{amsmath}
\usepackage{amsfonts}
\usepackage[latin1]{inputenc}
 \usepackage{graphicx}
    \usepackage{caption}
    \usepackage[skip=0.5ex, belowskip=1ex,
                labelformat=brace,
                singlelinecheck=off]{subcaption}
                
\usepackage{rotating}
\usepackage{float}
\usepackage{color}
\usepackage{tikz}
\usetikzlibrary{intersections}
\usepackage{keycommand}
  \usetikzlibrary{intersections,automata,arrows,positioning,calc}
\usepackage{pgfplots}
\usepackage{hyperref}

\newcommand{\BigI}{{\rm 1\kern-.17em 1}}
\newcommand{\BigZ}{{\rm Z\!\!Z}}
\newcommand{\BigR}{{\rm I\kern-.17em R}}
\newcommand{\BigN}{{\rm I\kern-.17em N}}
\newcommand{\BigC}{{\rm \kern.24em \vrule width.02em height1.4ex depth-.05ex \kern-.26em C}}

\newcommand{\R}{\mathbb  R}

\newcommand{\1}{\mathbb  I}

\newtheorem{lemma}{Lemma}[section]

\newtheorem{remark}[lemma]{Remark}

\newtheorem{definition}[lemma]{Definition}

\def\beq{\begin{equation}}   \def\eeq{\end{equation}}
\def\bea{\begin{eqnarray}}  \def\eea{\end{eqnarray}}

\newcommand\mysection{\setcounter{equation}{0}\section}
\renewcommand{\theequation}{\thesection.\arabic{equation}}
\newcounter{hran} \renewcommand{\thehran}{\thesection.\arabic{hran}}

\def\bmini{\setcounter{hran}{\value{equation}}
    \refstepcounter{hran}\setcounter{equation}{0}
    \renewcommand{\theequation}{\thehran\alph{equation}}\begin{eqnarray}}

\def\bminiG#1{\setcounter{hran}{\value{equation}}
\refstepcounter{hran}\setcounter{equation}{-1}
\renewcommand{\theequation}{\thehran\alph{equation}}
\refstepcounter{equation}\label{#1}\begin{eqnarray}}
 
\begin{document}

\title {  Numerical approaches to compute spectra of non-self adjoint operators in dimensions two and three}
\author{ Fatima Aboud \\
Mathematics Department, College of Science, University of Diyala, Iraq \\
{\it \small Fatima.Aboud@uodiyala.edu.iq} \\
Fran\c{c}ois Jauberteau \&  Didier Robert \\
Laboratoire de Math\'ematiques Jean Leray, CNRS-UMR 6629, \\
 Nantes Universit\'e, France \\
{\it \small francois.jauberteau@univ-nantes.fr,\; 
 didier.robert@univ-nantes.fr}
}
\vskip 1 truecm
\date{}
\maketitle
%
%
%

\begin{abstract}
In this article we are interested for the numerical computation of spectra of non-self adjoint quadratic operators, in two and three spatial dimensions. Indeed, in the multidimensional case very few results are known on the location of the eignevalues. This leads to solve nonlinear eigenvalue problems. In introduction we begin with a review of theoretical results and numerical results obtained for the one dimensional case. Then we present the numerical methods developed to compute the spectra (finite difference discretization) for the two and three dimensional cases. The numerical results obtained are presented and analyzed. One difficulty here is that we have to compute eigenvalues of strongly non-self-adjoint operators which are unstable. This work is in continuity of a previous work in one spatial dimension \cite{AFRM2021}.
\end{abstract}

\noindent {\bf Keywords~:} non-self adjoint quadratic operators, nonlinear eigenvalue problems, spectra, finite difference methods.

\tableofcontents

\pagestyle{myheadings}

\mysection{Introduction}

The aim of this paper is the extension of \cite{AFRM2021} (one dimensional in space) to higher dimensional spaces. \\

In this work we study the polynomial family of operators $L(\lambda) = H_{0}+ \lambda H_{1}+\cdots+
\lambda^{m-1}H_{m-1} + \lambda^{m} \1$, where the coefficients $H_{0}, H_{1}, \cdots, H_{m-1}$ are operators defined on some Hilbert space $\mathcal{H}$, $\1$ is the identity operator and $\lambda$ is a complex parameter. We are interested to study
the spectrum of the family $L(\lambda)$.\\

The problem $L(\lambda)u(x) = 0$, is called a non-linear eigenvalue problem for $m \geq 2$.  
The number $\lambda_{0} \in \mathcal{C}$ is called an eigenvalue of $L(\lambda)$, if there exists  $u_{0} \in\mathcal{H}$, $u_{0}\neq 0$ such that $L(\lambda_{0}) u_{0}=0$.
We consider here a quadratic family ($m = 2$) and in particular we are interested
in the case $L_{P} (\lambda)= \Delta +(P(x) \lambda)^{2}$, which is defined on the Hilbert space $\mathcal{L}^{2}(\mathbb{R}^{n})$, where $P$ is an elliptic positive polynomial of degree $M \geq 2$. For this example the known results for the existence of eigenvalues are given for $n = 1$ and $n$ is even (see \cite{herowa}) .\\

The main goal of the work of \cite{ab} and \cite{abro} is to check the following conjecture, stated by
Helffer-Robert-Wang \cite{herowa} :
\textit{for every dimension $n$, for every $M \geq 2$, the spectrum of $L_{P}$ is non empty}, which is confirmed for the following cases :
\begin{itemize}
\item $n = 1, 3$, for every polynomial $P$ of degree $M \geq 2$;
\item $n = 5$, for every convex polynomial $P$ satisfying some technical conditions;
\item $n = 7$, for every convex polynomial $P$.
\end{itemize}
This result extends to the case of quasi-homogeneous polynomial and quasi-elliptic, for example $P(x, y) = x^{2} + y^{4}$, $x \in \mathbb{R}^{n_{1}}$, $y \in \mathbb{R}^{n_{2}}$, $n_{1}+n_{2}=n$, and $n$ is
even. These last results were proved by computing the coefficients of a semi-classical trace formula and by using the theorem of Lidskii.\\

In \cite{AFRM2021}, some theoretical results are presented for $n=1$, for quadratic family of operators~:
$$
L(\lambda)=L_0+\lambda L_1 +\lambda^2 \1
$$
where $L_0$ and $L_1$ are operators on an Hilbert space and $\1$ is the identity operator. Moreover some numerical methods are presented to compute the spectrum of such operators with a  non self-adjoint linear eigenvalue problem, in one dimension in space. 
The proposed numerical methods are based on spectral methods and finite difference methods, in bounded and unbounded domains.
For bounded domain  homogeneous Dirichlet boundary conditions and periodic boundary conditions are considered with a comparison of the results obtained in unbounded and bounded domains. 
These comparisons show the difficulties for the numerical computation of the spectra of such operators. 
The numerical results obtained are in agreement with theoretical results on the location of the eigenvalues (see \cite{Christ}). \\

Here we are interested in the location in the complex plane of the eigenvalues for non-linear eigengenvalue problems in two and three dimensional cases. 
Indeed, in the multidimensional case very few theoretical results are known on the location of the eigenvalues for non-self adjoint quadratic operators. So numerical methods are developed to compute the spectra and numerical results are presented. \\

Here we notice some references of papers concerning computations of eigenvalues, connected  with our paper. For estimation of the spectra of linear differential operators we can refer to \cite{BrSch}, \cite{BrOs}, \cite{Osb}, \cite{ChaLem}, \cite{KO}, \cite{ALAM}, \cite{HANSEN} and \cite{Cha}. For numerical methods to solve nonlinear eigenvalue problems (finite dimensional case) and for linearization of such nonlinear eigenvalue problems we can refer to \cite{asle}, \cite{NLE1}, \cite{NLE2}, \cite{NLE3},\cite{NLE4}, \cite{NLE5}, \cite{NLE6} and \cite{NLE7}.\\

%
%
%

%
%
\newpage
\include{section1}

\include{section2}



\include{Numericalresults}

\mysection{Conclusion}

In this work we have considered quadratic family of operators~:
$$
L(\lambda)=L_0+\lambda L_1 +\lambda^2 \1
$$
where $L_0$ and $L_1$ are operators on an Hilbert space and $\1$ is the identity operator. \\

We have presented numerical methods to compute the spectrum of such operators. We reduce it to a  non self-adjoint linear eigenvalue problem. The numerical methods proposed are based on finite difference methods, in bounded domains. We consider homogeneous Dirichlet boundary conditions and periodic boundary conditions. \\

The numerical results obtained in the two and three dimensional cases are presented since very few theoretical results are known on the location of the eigenvalues in the multidimensional case. Numerical results obtained here show that theoretical results on the location of the eigenvalues in one dimension in space are still satisfied in dimensions 2 and 3. \\

\noindent{\bf Acknowledgments} \\

\noindent This work has benefit of the visit of Fatima Aboud, at the Laboratoire de Math\'ematiques Jean Leray,
Nantes Universit\'e (France). This visit was supported by the F\'ed\'eration de Math\'ematiques des Pays de la Loire, CNRS FR 2962. \\

%
%

%
\include{figures}

\end{document}

%% file: section1.tex
\section{Homogeneous Dirichlet boundary conditions}
\subsection{The two dimensional case}
We consider $n=2$ and the operator
$$L_\lambda(u)=-\Delta u+(|x|^2-\lambda)^2u$$
where 
$$\Delta u(x,y)=\frac{\partial^2u}{\partial {x^2}} +\frac{\partial^2u}{\partial {y^2}}, \; (x,y)\in \R^2.$$
We discretize the problem by using finite difference method to solve the following homogeneous problem with the homogeneous Dirichlet boundary conditions
\begin{equation}\label{HHDC2}
\begin{cases}
L_\lambda(u(x,y))=0,   & \text{$(x,y) \in [-L,L]\times [-L,L]$} \\
u(L,y)=u(-L,y)=0  & \text{ }\\
u(x,L)=u(x,-L)=0  & \text{ }\\
\end{cases}
\end{equation}
So we take the rectangular grid $[-L,L]\times [-L,L]$ with $\displaystyle h=\frac{2L}{N-1}$ for some $N$ ($N$ the number of points on the interval $[-L,L]$). To approximate the Laplacian we need to $5-$points $(i,j+1)$, $(i,j-1)$, $(i-1,j)$, $(i+1,j)$ and $(i,j)$ (see Figure \ref{5-points}). For simplicity we denote the point $(x_i,y_j)$ by $(i,j)$ and $u(x_i,y_j)$ by $u_{i,j}$.\\

%
            We approximate the derivatives $\displaystyle \frac{\partial^2u}{\partial x^2}$, $\displaystyle \frac{\partial^2u}{\partial y^2 }$ by using Taylor expansion so~:
$$\frac{\partial^2u}{\partial {x^2}}|_{i}= \frac{u_{i-1,j}-2u_{i,j} +u_{i+1,j}}{h^2}+O(h^2)$$
$$ \frac{\partial^2u}{\partial {y^2}}|_{j} = \frac{u_{i,j-1}-2u_{i,j}+u_{i,j+1}}{h^2}+O(h^2) $$
we get
$$\Delta u(x,y)|_{i,j}\approx \frac{u_{i-1,j}-2u_{i,j} +u_{i+1,j}}{h^2}+\frac{u_{i,j-1}-2u_{i,j} +u_{i,j+1}}{h^2} $$
$$\Delta u(x,y)|_{i,j}\approx \frac{-4u_{i,j}+u_{i-1,j} +u_{i+1,j}+u_{i,j-1}+u_{i,j+1}}{h^2}. $$

Now we define 
  $$ U_{\bullet j}=\left(\begin{array}{c}
                  u_{1,j}\\
                  u_{2,j}\\
                  \vdots\\
                  u_{N,j}
                   \end{array}\right)_{N\times 1}
, \;\;\;
 \overrightarrow{U}  =\left(\begin{array}{c}
                  U_{\bullet 1}\\
                 U_{\bullet 2} \\
                  \vdots\\
                 U_{\bullet N} 
                   \end{array}\right)_{N^2\times 1}
$$
$$ X_{\bullet j}=\left(\begin{array}{c}
                  {x_1,y_j}\\
                  {x_2,y_j}\\
                  \vdots\\
                  {x_N,y_j}
                   \end{array}\right)_{N\times 1}
, \;\;\;
 \overrightarrow{X}  =\left(\begin{array}{c}
                  X_{\bullet 1}\\
                 X_{\bullet 2} \\
                  \vdots\\
                 X_{\bullet N} 
                   \end{array}\right)_{N^2\times 1}
$$
$$ |X_{\bullet j}|^2=\left(\begin{array}{c}
                  {x_1^2+y_j^2}\\
                  {x_2^2+y_j^2}\\
                  \vdots\\
                  {x_N^2+y_j^2}
                   \end{array}\right)_{N\times 1}
, \;\;\;
 |\overrightarrow{X}|^2  =\left(\begin{array}{c}
                 | X_{\bullet 1}|^2\\
                | X_{\bullet 2}|^2 \\
                  \vdots\\
                | X_{\bullet N}| ^2
                   \end{array}\right) _{N^2\times 1}
$$         

To construct the Laplacian matrix we need to define the matrix $\mathbb{D}_4$ which is an $N\times N$ matrix such that  
$$\mathbb{D}_4 = \left( \begin{array}{ccccccc}

         -4 & 1    &0 &\cdots  &  & \cdots & 0 \\

          1 & -4   &1 & 0 &\cdots  & \cdots & 0 \\
0 & 1 &-4&1    &  0  & \cdots & 0 \\
 \vdots  &   & \vdots  &  & \vdots &  &   \vdots  \\
 0 & \cdots &  0& 1 &- 4  &    1 &0 \\
0 & \cdots & \cdots &0  &1  &    - 4 &1 \\
          0 & \cdots &  &\cdots  & 0  &        1  &-4 \\
\end{array} \right)_{N \times N}$$
We denote the Laplacian matrix by $\underline{\Delta}$ which has the following form
\begin{equation}\label{Laplacian2D}
h^2 \underline{\Delta} = \left( \begin{array}{ccccccc}

          \mathbb{D}_4  & \1    &0 &\cdots  &  & \cdots & 0 \\

          \1 & \mathbb{D}_4   &\1 & 0 &\cdots  & \cdots & 0 \\
0 & \1 &\mathbb{D}_4&\1    &  0  & \cdots & 0 \\
 \vdots  &   & \vdots  &  & \vdots &  &   \vdots  \\
 0 & \cdots & 0 & \1 & \mathbb{D} _4 &     \1 & 0 \\
0 & \cdots & \cdots &0  & \1  &     \mathbb{D}_4 & \1 \\
          0 & \cdots &  &\cdots  & 0  &        \1  &\mathbb{D}_4\\
\end{array} \right)_{N^2 \times N^2}
\end{equation}
where $\1$  is the identity $N\times N$ matrix and we note that the zeros represent the zero $N\times N$ matrices. So we have the following discrete problem
\begin{equation}
 -\underline{\Delta}  \overrightarrow{U}+(|\overrightarrow{X}|^2-\lambda \1)^2\overrightarrow{U} =0
\end{equation}
i.e.
$$
 \Big(-\underline{\Delta}+ diag(|\overrightarrow{X}|^4)\Big)  \overrightarrow{U}-2\lambda\; diag(|\overrightarrow{X}|^2) \;  \;\overrightarrow{U}+\lambda^2\; \1 \; \overrightarrow{U} =0
$$
We set 
$$ {H}_0= -\underline{\Delta}+\; \; diag(|\overrightarrow{X}|^4) , \;\;\; {H}_1= -2 \; diag( |\overrightarrow{X}|^2)$$
where $diag(|\overrightarrow{X}|^4)$, $diag( |\overrightarrow{X}|^2)$ represent  square diagonal matrices with the elements of vector $|\overrightarrow{X}|^4$, $ |\overrightarrow{X}|^2$ on the main diagonal, respectively.
By using the linearization of this problem we get 
\begin{equation}
 \mathcal{A}= \left( \begin{array}{cc}
                        0   & \1  \\
                      -H_0  & - H_1
                     \end{array} \right)
\end{equation}

\subsection{The three  dimensional case}
We consider the case of the space dimension $n=3$ for the family of operators:
$$L_\lambda=-\Delta+(|{\bf x}|^2-\lambda)^2 \1 =-\Delta+|{\bf x}|^4 \1- 2\lambda |{\bf x}|^2 \1+\lambda^2 \1$$
with ${\bf x} = (x,y,z) \in \BigR^3$ and $\1$ the identity operator. Introducing a real parameter $c$ we study a more general case of this family of operators
$$L_{\lambda,c}=-\Delta+|{\bf x}|^4-2c\lambda |{\bf x}|^2+\lambda^2$$
%
\begin{equation}\label{HHDC3}
\begin{cases}
L_{\lambda,c}(u(x,y,z))=0,   & \text{$(x,y,z) \in [-L,L]\times [-L,L]\times [-L,L]$} \\
u(L,y,z)=u(-L,y,z)=0  & \text{ }\\
u(x,L,z)=u(x,-L,z)=0  & \text{ }\\
u(x,y,L)=u(x,y,-L)=0  & \text{ }\\
\end{cases}
\end{equation}
%
For simplicity reason, we denote the point $(x_i,y_j, z_k)$ by $(i,j,k)$ and $u(x_i,y_j,z_k)$ by $u_{i,j,k}$.\\
To calculate the Laplacian we need to $7-$points $(i,j,k), (i-1,j,k), (i+1,j,k), (i,j-1,k), (i,j+1,k), (i,j,k-1), (i,j,k+1)$ (see Figure \ref{7-points}). We have~: 
$$\Delta u(x,y,z)=\frac{\partial^2u}{\partial {x^2}} +\frac{\partial^2u}{\partial {y^2}}+\frac{\partial^2u}{\partial {z^2}}, \; (x,y,z)\in \R^3$$
and we approximate the derivatives $\displaystyle \frac{\partial^2u}{\partial{x^2}}$, $\displaystyle \frac{\partial^2u}{\partial {y^2}}$ and $\displaystyle \frac{\partial^2u}{\partial{z^2}}$ by using Taylor expansion~:  
$$\frac{\partial^2u}{\partial {x^2}}|_{i}= \frac{u_{i-1,j,k}-2u_{i,j,k} +u_{i+1,j,k}}{h^2}+O(h^2)$$
$$ \frac{\partial^2u}{\partial {y^2}}|_{j} = \frac{u_{i,j-1,k}-2u_{i,j,k}+u_{i,j+1,k}}{h^2}+O(h^2) $$
$$ \frac{\partial^2u}{\partial {z^2}}|_{k} = \frac{u_{i,j,k-1}-2u_{i,j,k}+u_{i,j,k+1}}{h^2}+O(h^2) $$
So we obtain~:
$$\Delta u(x,y,z)|_{i,j,k}\approx \frac{u_{i-1,j,k}-2u_{i,j,k} +u_{i+1,j,k}}{h^2}+\frac{u_{i,j-1,k}-2u_{i,j,k}+u_{i,j+1,k}}{h^2}+ \frac{u_{i,j,k-1}-2u_{i,j,k}+u_{i,j,k+1}}{h^2} $$
{i.e.}
$$\Delta u(x,y,z)|_{i,j,k} \approx \frac{-6u_{i,j,k}+u_{i-1,j,k} +u_{i+1,j,k}+u_{i,j-1,k}+u_{i,j+1,k}+u_{i,j,k-1}+u_{i,j,k+1}}{h^2} $$
%
%
%
%
%
%

\begin{figure}[!tbp]
  \centering
  \begin{minipage}[b]{0.4\textwidth}

\begin{center}
\begin{tikzpicture}
\draw (-1.5,0) node[sloped,above]{$(i-1,j)$}circle (2pt) -- (1.5,0) node[sloped,above]{$(i+1,j)$}circle (2pt);
\draw (0,1.5) node[sloped,above]{$(i,j-1)$}circle (2pt) -- (0,0)  circle (2pt) -- (0,-1.5) node[sloped,below]{$(i,j+1)$}circle (2pt);

\draw (0.5,-0.65) node [above=0pt,fill=white]{$(i,j)$};
\end{tikzpicture}
\end{center}
    \caption{$5-$points finite difference $2D$}
\label{5-points}
  \end{minipage}
  \hfill
  \begin{minipage}[b]{0.4\textwidth}
   
   \begin{tikzpicture}
\draw (-3,0,0) node[sloped,above]{$(i-1,j,k)$}circle (2pt) -- (3,0,0) node[sloped,above]{$(i+1,j,k)$}circle (2pt);
\draw (0,3,0) node[sloped,above]{$(i,j-1,k)$}circle (2pt) -- (0,-3,0) node[sloped,below]{$(i,j+1,k)$}circle (2pt);
\draw (0,0,3)node[sloped,below]{$(i,j,k-1)$}circle (2pt)  -- (0,0,0)  circle (2pt) -- (0,0,-3) node[sloped,above]{$(i,j,k+1)$} circle (2pt);
\draw (0.75,-0.75,0) node [above=0pt,fill=white]{$(i,j,k)$};
\end{tikzpicture}
   
    \caption{$7-$points finite difference $3D$}
\label{7-points}
  \end{minipage}
\end{figure}
Now we define 
  $$ U_{\bullet, j,k}=\left(\begin{array}{c}
                  u_{1,j,k}\\
                  u_{2,j,k}\\
                  \vdots\\
                  u_{N,j,k}
                   \end{array}\right)_{N\times 1}
,\;\;\;
V_{i,\bullet,k}=\left(\begin{array}{c}
                  U_{\bullet,1,k}\\
                  U_{\bullet,2,k}\\
                  \vdots\\
                  U_{\bullet,N,k}
                   \end{array}\right)_{N^2\times 1}
                   ,\;\;\;
\overrightarrow{U}=\left(\begin{array}{c}
                  V_{i,\bullet,1}\\
                 V_{i,\bullet,2}\\
                  \vdots\\
                 V_{i,\bullet,N}
                   \end{array}\right)_{N^3\times 1}
$$
$$ X_{\bullet, j,k}=\left(\begin{array}{c}
                  {x_1,y_j,z_k}\\
                  {x_2,y_j,z_k}\\
                  \vdots\\
                  {x_N,y_j,z_k}
                   \end{array}\right)_{N\times 1}
, \;\;\;
X'_{i,\bullet,k}=\left(\begin{array}{c}
                  X_{\bullet, 1,k}\\
                  X_{\bullet, 2,k}\\
                  \vdots\\
                  X_{\bullet, 3,k}
                   \end{array}\right)_{N^2\times 1}
, \;\;\;
\overrightarrow{X}=\left(\begin{array}{c}
                  X'_{i,\bullet,1}\\
                  X'_{i,\bullet,2}\\
                  \vdots\\
                  X'_{i,\bullet,N}
                   \end{array}\right)_{N^3\times 1}
$$
$$
|\overrightarrow{X}|^2=\left(\begin{array}{c}
                  x_i^2+y_\bullet^2+z_1^2\\
                  x_i^2+y_\bullet^2+z_2^2\\
                  \vdots\\
                  x_i^2+y_\bullet^2+z_N^2
                   \end{array}\right)_{N^3\times 1}, \;\;\;\;
|\overrightarrow{X}|^4=\left(\begin{array}{c}
                 \Big( x_i^2+y_\bullet^2+z_1^2 \Big)^2\\
                 \Big( x_i^2+y_\bullet^2+z_2^2\Big)^2\\
                  \vdots\\
                 \Big( x_i^2+y_\bullet^2+z_N^2\Big)^2
                   \end{array}\right)_{N^3\times 1}                   
$$
To construct the matrix of Laplacian $\underline{\Delta}$ we need to define the following matrices. 
Let $\mathbb{D}_6$ be an ${N\times N}$ matrix defined as follows
$$\mathbb{D}_6 = \left( \begin{array}{ccccccc}

          -6 & 1    &0 &\cdots  &  & \cdots & 0 \\

          1 & -6   &1 & 0 &\cdots  & \cdots & 0 \\
0 & 1 &-6&1    &  0  & \cdots & 0 \\
 \vdots  &   & \vdots  &  & \vdots &  &   \vdots  \\
 0 & \cdots & 0 &1  & -6 &     1 &0 \\
0 & \cdots & \cdots &0  &1  &     -6 &1 \\
          0 & \cdots &  &\cdots  & 0  &        1  &-6 \\
\end{array} \right)_{N\times N}$$
Now we define $\mathbb{G}_6$ which is an ${N^2}\times {N^2}$ matrix    such that
$$\mathbb{G}_6 = \left( \begin{array}{ccccccc}

          \mathbb{D}_6 & \1    &0 &\cdots  &  & \cdots & 0 \\

          \1 & \mathbb{D}_6   &\1 & 0 &\cdots  & \cdots & 0 \\
0 & \1 &\mathbb{D}_6 &\1    &  0  & \cdots & 0 \\
 \vdots  &   & \vdots  &  & \vdots &  &   \vdots  \\
 0 & \cdots &  0& \1  &\mathbb{D}_6 &     \1  &0 \\
0 & \cdots & \cdots &0  &\1  &     \mathbb{D}_6  &\1 \\
          0 & \cdots &  &\cdots  & 0  &        \1  &\mathbb{D}_6  \\
\end{array} \right)_{N^2\times N^2}$$
where $ \1$ is the $N\times N$ identity and the zeros represent the $N\times N$ zero matrice.
So the $3D$ Laplacian finite difference matrix is given by
\begin{equation}\label{Lap3D}
h^2 \underline{\Delta }= \left( \begin{array}{ccccccc}

          \mathbb{G}_6 & \mathbb{J}   &0 &\cdots  &  & \cdots & 0 \\

          \mathbb{J} & \mathbb{G}_6  &\mathbb{J} & 0 &\cdots  & \cdots & 0 \\
0 & \mathbb{J}&\mathbb{G}_6&\mathbb{J}   &  0  & \cdots & 0 \\
 \vdots  &   & \vdots  &  & \vdots &  &   \vdots  \\
 0 & \cdots & 0 & \mathbb{J} &\mathbb{G}_6  &     \mathbb{J}  &0\\
0 & \cdots & \cdots &0  &\mathbb{J}  &     \mathbb{G}_6 &\mathbb{J} \\
          0 & \cdots &  &\cdots  & 0  &        \mathbb{J}  &\mathbb{G}_6 \\
\end{array} \right)_{N^3\times N^3}
\end{equation}
where $\mathbb{J}$ is the $N^2\times N^2$ identity and the zeros are the zero $N^2\times N^2$ matrice. We note that the Laplacian matrix will be $N^3\times N^3$.\\

So we have the following discrete problem
\begin{equation}
 -\underline{\Delta}  \overrightarrow{U}+(|\overrightarrow{X}|^2-\lambda \1)^2\overrightarrow{U} =0
\end{equation}
With the parameter $c$ we deduce the more general discrete problem 
$$
 \Big(-\underline{\Delta}+diag(|\overrightarrow{X}|^4)\Big)  \overrightarrow{U}-2\lambda\;c \; diag(|\overrightarrow{X}|^2) \;  \;\overrightarrow{U}+\lambda^2\; \1 \; \overrightarrow{U} =0
$$
We set 
$$ {H}_0= -\underline{\Delta}+diag(|\overrightarrow{X}|^4), \;\;\; {H}_1= -2\;c\;diag(|\overrightarrow{X}|^2)$$
By using the linearization of this problem we get 
\begin{equation}
 \mathcal{A}= \left( \begin{array}{cc}
                        0   & \1  \\
                      -H_0  & - H_1
                     \end{array} \right)
\end{equation}
with  
$$H_0 = \left( \begin{array}{cccccc}

        \Theta_{i,\bullet,1}- \mathbb{G}_6/h^2  & -\mathbb{J}   & 0  & \cdots & \cdots & 0 \\
\\
          -\mathbb{J} &  \Theta_{i,\bullet,2}-\mathbb{G}_6/h^2   & -\mathbb{J} & 0  & \cdots & 0 \\
          \\
 \vdots  & \vdots  & \vdots & \vdots &  \vdots &   \vdots  \\
 \\
0 & \cdots & 0  &-\mathbb{J}  &    \Theta_{i,\bullet,N-1}-\mathbb{G}_6/h^2 &-\mathbb{J} \\
\\
          0 & \cdots &  \cdots  & 0  &      - \mathbb{J}  &  \Theta_{i,\bullet,N}-\mathbb{G}_6/h^2 \\
\end{array} \right)_{N^3\times N^3}$$
where $\Theta_{i,\bullet,j}=diag(|X_{i,\bullet,j}|^4)$ is the diagonal $N^2\times N^2$ matrix with $|X'_{i,\bullet,j}|^4$ on the diagonal, and
$H_1=-2 c \;\; diag(|\overrightarrow{X}|^2)$  
with $diag(|\overrightarrow{X}|^2)$ is the diagonal $N^3\times N^3$ matrix with $|\overrightarrow{X}|^2$ on the diagonal.\\

%
%

%% file: section2.tex
\section{Periodic case}
Consider the following differential operator~:
\begin{equation}\label{periodicOp}
 -\Delta u + \lambda^2 u+ \frac{\lambda}{i} \sum_{j=1}^n a_j(x_j) \frac{\partial u}{\partial x_j}=0, \;\;x\in \mathbb{R}^n
\end{equation}
with periodic boundary conditions.\\

In the case where the coefficients $a_j(x_j)$ are space independent, 
if we consider $u=\sum _{{\bf k}\in \BigZ}\hat{u}_{\bf k} e^{i{\bf k} \cdot {\bf x} }$ and substituting in (\ref{periodicOp}) we obtain~:
$$k^2\hat{u}_{\bf k}+\lambda^2 \hat{u}_{\bf k} +\frac{\lambda}{i} \sum_{j=1}^n i a_j k_j \hat{u}_{\bf k}=0$$
$$\Rightarrow \lambda^2 +\lambda \sum_{j=1}^n a_j k_j + k^2 =0 $$
where ${\bf x} = (x_i)_{i=1,\ldots,n}$, ${\bf k} = (k_i)_{i=1,\ldots,n}$ and $k = \vert {\bf k}\vert$. \\

In the following we consider the case where the coefficients are space dependent and periodic in space. 

\subsection{The two dimensional case}
Let $\Omega=(0,L)\times(0,L)\subset\mathbb{R}^2$, we study the following problem
\begin{eqnarray*}\label{periodicOp2D}
 -\Delta u+ \lambda^2 u+ \frac{\lambda}{i} \Big( a_1(x) \frac{\partial u}{\partial x} + a_2(y) \frac{\partial u}{\partial y}\Big)=0\\
 u(x+L,y+L)=u(x,y)\;\;\;\\  
\end{eqnarray*}
%
%
where $L\in \R$, $x\in [0,L]$ and $ y\in [0,L]$. So we take the rectangular grid $[0,L]\times [0,L]$ with $\displaystyle h=\frac{2L}{N-1}$ for some $N$ ($N$ the number of points on the interval $[0,L]$). 
We set~:
$$H_0=-\Delta  \; ,\;\; H_1=\frac{1}{i} \sum_{j=1}^2 a_j(x_j) \frac{\partial }{\partial x_j}$$
To get the formula of the linearization operator $\mathcal{A}$ we need to have the formulas of the matrices $H_0, H_1$. We note that
$H_0$ is the Laplacian, so it has the formula given in (\ref{Laplacian2D}). 
To calculate the Jaccobi matrix we need to $4-$points $(i,j+1)$, $(i,j-1)$, $(i-1,j)$, $(i+1,j)$ (for simplicity we denote the point $(x_i,y_j)$ by $(i,j)$ and $u(x_i,y_j)$ by $u_{i,j}$).
We approximate the derivatives $\displaystyle \frac{\partial u}{\partial  x }$, $\displaystyle \frac{\partial u}{\partial  y }$ by using Taylor expansion so 
$$\frac{\partial u}{\partial {x}}|_{i}= \frac{u_{i+1,j}-u_{i-1,j}}{2h}+O(h^2)$$
$$ \frac{\partial u}{\partial {y}}|_{j} = \frac{u_{i,j+1}-u_{i,j-1}}{2h}+O(h^2) $$
%
 %
%
%
%
%
To give the formula of $H_1$ we need the following matrices
$$\mathbb{D}_x = \left( \begin{array}{ccccccc}
       0 & a_1(x_1)    &0 &\cdots  &  & \cdots & -a_1(x_1) \\
          -a_1(x_2)  & 0  & a_1(x_2) & 0 &\cdots  & \cdots & 0 \\
0 & -a_1(x_3)  &0&a_1(x_3)     &  0  & \cdots & 0 \\
 \vdots  &   & \vdots  &  & \vdots &  &   \vdots  \\
 0 & \cdots &  0& -a_1(x_{N-2})  &0  &     a_1(x_{N-2}) &0 \\
0 & \cdots & \cdots &0  & -a_1(x_{N-1})  &    0 & a_1(x_{N-1}) \\
          a_1(x_N) & \cdots &  &\cdots  & 0  &       -a_1(x_{N})    &0 \\
\end{array} \right)_{N \times N}$$
so
\begin{equation}\label{}
H_1 = \frac{1}{2 i h} \left( \begin{array}{ccccccc}

          \mathbb{D}_x &  {Dy_1}   &0 &\cdots  &  & \cdots & {-Dy_1} \\

         - {Dy_2}  & \mathbb{D}_x   &{Dy_2}  & 0 &\cdots  & \cdots & 0 \\
0 & -{Dy_3}  &\mathbb{D}_x&{Dy_3}    &  0  & \cdots & 0 \\
 \vdots  &   & \vdots  &  & \vdots &  &   \vdots  \\
 0 & \cdots & 0 &  -{Dy_{N-2}} & \mathbb{D} _x &   {Dy_{N-2}} & 0 \\
0 & \cdots & \cdots &0  &  -{Dy_{N-1}} &     \mathbb{D}_x &  {Dy_{N-1}} \\
          {Dy_{N}} & \cdots &  &\cdots  & 0  &       - {Dy_{N}}  &\mathbb{D}_x\\
\end{array} \right)_{N^2 \times N^2}
\end{equation}
where ${Dy_j}= diag(a_2(y_{j}))$.
\subsection{The three dimensional case}
Let $\Omega=(0,L)\times(0,L)\times(0,L)\subset\mathbb{R}^3$, we study the following problem
\begin{eqnarray*}\label{periodicOp3D}
 -\Delta u+ \lambda^2 u+ \frac{\lambda}{i} \Big( a_1(x) \frac{\partial u}{\partial x} + a_2(y) \frac{\partial u}{\partial y}+ a_3(z) \frac{\partial u}{\partial z}\Big)=0\\
 \;\;\;u(x+L,y+L,z+L)=u(x,y,z)\hfill\\ 
\end{eqnarray*}
%
%
where $x\in [0,L]$, $ y\in [0,L]$, $ z\in [0,L]$ and $L\in \R$. So we take the rectangular grid $[0,L]\times [0,L]\times [0,L]$ with $\displaystyle h=\frac{2L}{N-1}$ for some $N$ ($N$ the number of points on the interval $[0,L]$). We set~: 
$$H_0=-\Delta \; ,\;\; H_1=\frac{1}{i} \sum_{j=1}^3 a_j(x_j) \frac{\partial }{\partial x_j}$$
To get the formula of the linearization operator $\mathcal{A}$ we need to have the formulas of the matrices $H_0, H_1$. We note that 
$H_0$ is the Laplacian, so it has the formula given in (\ref{Lap3D}). \\
For simplicity reason, we denote the point $(x_i,y_j, z_k)$ by $(i,j,k)$ and $u(x_i,y_j,z_k)$ by $u_{i,j,k}$.\\
To calculate the Jaccobi matrix we need to $6-$points $(i-1,j,k), (i+1,j,k), (i,j-1,k), (i,j+1,k), (i,j,k-1), (i,j,k+1)$. 
We approximate the derivatives $\displaystyle \frac{\partial u}{\partial x }$, $\displaystyle \frac{\partial u}{\partial y }$ and $\displaystyle \frac{\partial u}{\partial z}$ by using Taylor expansion so~:
$$\frac{\partial u}{\partial {x}}|_{i}= \frac{u_{i+1,j,k}-u_{i-1,j,k}}{2h}+O(h^2)$$
$$ \frac{\partial u}{\partial {y}}|_{j} = \frac{u_{i,j+1,k}-u_{i,j-1,k}}{2h}+O(h^2) $$
$$ \frac{\partial u}{\partial {z}}|_{k} = \frac{u_{i,j, k+1}-u_{i,j, k-1}}{2h}+O(h^2) $$
%
 %
%
%
%
%
To give the formula of $H_1$, we need the following matrices
$$\mathbb{D}_x = \left( \begin{array}{ccccccc}

         0 & a_1(x_1)    &0 &\cdots  &  & \cdots & -a_1(x_1) \\

          -a_1(x_2)  & 0  & a_1(x_2) & 0 &\cdots  & \cdots & 0 \\
0 & -a_1(x_3)  &0&a_1(x_3)     &  0  & \cdots & 0 \\
 \vdots  &   & \vdots  &  & \vdots &  &   \vdots  \\
 0 & \cdots &  0& -a_1(x_{N-2})  &0  &     a_1(x_{N-2}) &0 \\
0 & \cdots & \cdots &0  & -a_1(x_{N-1})  &    0 & a_1(x_{N-1}) \\
          a_1(x_N) & \cdots &  &\cdots  & 0  &       -a_1(x_{N})    &0 \\
\end{array} \right)_{N \times N}$$
$${Dy_j}= \Big(diag(a_2(y_{j}))\Big)_{N \times N}\; \mbox{and} \;\; {Dz_j}= \Big(diag(a_3(z_{j}))\Big)_{N^2 \times N^2}$$
so 
\begin{equation}
\mathbb{G} _{x,y} = \left( \begin{array}{ccccccc}

          \mathbb{D}_x &  {Dy_1}   &0 &\cdots  &  & \cdots & -{Dy_1} \\

         - {Dy_2}  & \mathbb{D}_x   &{Dy_2}  & 0 &\cdots  & \cdots & 0 \\
0 & -{Dy_3}  &\mathbb{D}_x&{Dy_3}    &  0  & \cdots & 0 \\
 \vdots  &   & \vdots  &  & \vdots &  &   \vdots  \\
 0 & \cdots & 0 &  -{Dy_{N-2}} & \mathbb{D} _x &  {Dy_{N-2}} & 0 \\
0 & \cdots & \cdots &0  &  -{Dy_{N-1}} &     \mathbb{D}_x &  {Dy_{N-1}} \\
          {Dy_{N}} & \cdots &  &\cdots  & 0  &       - {Dy_{N}}  &\mathbb{D}_x\\
\end{array} \right)_{N^2 \times N^2}
\end{equation}
thus
\begin{equation}
H_1 = \frac{1}{2 i h}\left( \begin{array}{ccccccc}

          \mathbb{G} _{x,y} &  {Dz_1}   &0 &\cdots  &  & \cdots & - {Dz_1} \\

         - {Dz_2}  & \mathbb{G} _{x,y}  &{Dz_2}  & 0 &\cdots  & \cdots & 0 \\
0 & -{Dz_3}  &\mathbb{G} _{x,y}&{Dz_3}    &  0  & \cdots & 0 \\
 \vdots  &   & \vdots  &  & \vdots &  &   \vdots  \\
 0 & \cdots & 0 &  -{Dz_{N-2}} & \mathbb{G} _{x,y} & {Dz_{N-2}} & 0 \\
0 & \cdots & \cdots &0  &  -{Dz_{N-1}} &     \mathbb{G} _{x,y} & {Dz_{N-1}} \\
          {Dz_{N}} & \cdots &  &\cdots  & 0  &       - {Dz_{N}}  &\mathbb{G} _{x,y}\\
\end{array} \right)_{N^3 \times N^3}
\end{equation}

%% file: Numericalresults.tex
\section{Numerical results and discussion}

In this section we present the numerical results obtained in dimensions 2 and 3 for homogeneous boundary conditions and for periodic boundary conditions. The numerical computation of the spectra of non-self adjoint operators being unstable, it is possible to compute the pseudospectra, see for example \cite{Davies}, \cite{hiti1}, \cite{hiti2}, \cite{thr1}, \cite{trem}. Indeed the computation of the pseudospectra is more stable, but it is expensive since it requires to compute the minimal singular value $s_{min} (\mathcal{A}-z \1)$ at each point $z$ of the mesh. Here we have computed the spectra using Lapack library.

\subsection{ Homogeneous boundary conditions }

In this section we present the numerical results obtained for homogeneous boundary conditions and spatial dimensions $n=2$ and $n=3$. The operator considered here is $L_{\lambda,c}$ defined in (\ref{HHDC3}).\\

We start by presenting the results for $c=0$ and $N=100$ (resp. $N=25$) for the two (resp. three) dimensional case. In this case, the spectrum contains just a purely imaginary eigenvalues, see Figure \ref{coefc0dim23}. \\

When the parameter $c$ is increased ($c=0.5$, $1$, $2$, $5$ and $10$) the real part of the eigenvalues $\lambda$ is increased, see Figures \ref{coefc1dim23} \textendash \ref{coefc20dim23}. In one spacial dimension and $c=1$ according the theory we have that $Sp(L_{\lambda,c})$ is included in the two sectors $\{ \lambda \in \BigC, \vert arg(\lambda)\vert \ge \frac{\pi}{3} \}$ (see \cite{Christ} and \cite{AFRM2021} for a review of theoretical results).  Numerical results obtained here show that this location of the eigenvalues is still satisfied in dimensions 2 and 3. \\

The CPU time is expensive in 2D and in 3D, since the size of the matrices obtained after linearization and discretization is $2N^n\times 2N^n$, with $n=2$ and $n=3$. So the spatial resolution is quite limited. We have studied the influence of the spatial resolution on the mumerical results for $c=1$. On Figures \ref{compare2Dsamediff} and \ref{compare3Dsamediff} we test the influence of the choice of the parameter $N$ with a constant size of the domain $[-L,L]$. We can observe a similar comportment in $2D$ and in $3D$. When $N$ is increased the repartition of the imaginary part of the eigenvalues is increased but it is not the case for the real part. However on Figures \ref{compar2Ddiffsam} and \ref{compar3Ddiffsam} when the resolution $N$ is constant and the size of the domain is increased ($L=0.5$, $1$ and $2$) we can see that the repartition of the real and imaginary parts of the eigenvalues are increased. So, in conclusion increasing the parameter $N$ has an effect on the repartition of the imaginary part. Moreover, increasing the size of the domain has an effect on the repartition on the real and imaginary parts of the eigenvalues. \\

In conclusion, increasing $c$ increased the repartition of the spectrum along the real part. Increasing $N$ increased the repartition of the spectrum along the imaginary part. And increasing $L$ increased the repartition of the spectrum along the real and imaginary parts. \\




\subsection{ Periodic boundary conditions }

In this section we present the numerical results obtained for periodic boundary conditions and spatial dimensions $n=2$ and $n=3$. The operator considered here is defined in (\ref{periodicOp}). The domain retained for the simulations is $[-\pi,\pi]$. For the spatial resolution parameter $N$ we have retained $N=100$ (resp. $N=25$) for the simulations in one (resp. two) dimension in space. And we consider different choices for the coefficients $a_j(x_j)$, $j=1,\ldots,n$, periodic in space. \\

Firstly we study the two dimensional case with constant coefficients $a_j(x_j)$, $j=1,2$. We have presented the results obtained with different choices of the coefficients~: $a_1(x_1)=a_2(x_2)=1$ on Figure \ref{dim2periodic11}, $a_1(x_1)=1$, $a_2(x_2)=\sqrt{2}$ on Figure \ref{dim2periodic1sqrt2} and $a_1(x_1)=1$, $a_2(x_2)=5\sqrt{2}$ on Figure \ref{dim2periodic15sqrt2}. We can observe that the repartition of the real (resp. imaginary) parts of the eigenvalues are symetric with respect to the vertical (resp. horizontal) axes, which was not the case for the previous operator (\ref{HHDC3}). Moreover if the repartition of the eigenvalues is quite similar on Figures \ref{dim2periodic11} and \ref{dim2periodic1sqrt2} (the real part increases with the imaginary part in modulus), it is different for the eigenvalues with small imaginary part on Figure \ref{dim2periodic15sqrt2} (the real part decreases when the imaginary part increases in modulus). \\

Then we consider the two dimensional case with sinusoidal coefficients $a_1(x_1)=\sin(x_1)$ and $a_2(x_2)=\sin(x_2)$. The numerical results are presented on Figure \ref{dim2periodicsinsin}. We can see that the repartition of the eigenvalues is always symetric, as for constant coefficients. However the real part of the eigenvalues are smaller in modulus than for constant coefficients. Moreover the real part of the eigenvalues is approximatively equal to zero when the modulus of the imaginary part is less than $20$. \\

Now we present the three dimensional case with different choices of the coefficients. Numerical results for constant coefficients $a_j(x_j)=1$, $i=1,\ldots3$ and $a_1(x_1)=1$, $a_2(x_2)=\sqrt{2}$, $a_3(x_3)=5\sqrt{2}$ are represented on Figures \ref{dim3periodic111} and \ref{dim3periodic1sqrt25sqrt2} respectively. Then we consider sinusoidal coefficients $a_1(x_1)=\sin(x_1)$, $a_2(x_2)=\sin(x_2)$ and $a_3(x_3)=\sin(x_3)$. The numerical results are presented on Figure \ref{dim3periodicsinsinsin}. As for the periodic two dimensional case, we can see that the repartition of the eigenvalues is symetric with respect the horizontal and vertical axes when the coefficients $a_j(x_j)=1$, $i=1,\ldots3$ are the same (see Figures \ref{dim3periodic111} and \ref{dim3periodicsinsinsin}). However, when these coefficients are different these symetries are lost for the eigenvalues with small imaginary part (symetry is with respect the origin, see Figure \ref{dim3periodic1sqrt25sqrt2}). Moreover, apart from the eigenvalues with small imaginary parts, the eigenvalues are localized on the imaginary axis. 



%% file: figures.tex
\section{Figures}

\begin{figure}[hbt!]
\centering
 \begin{subfigure}{0.40\textwidth}
  \centering
  \includegraphics[width=\linewidth]{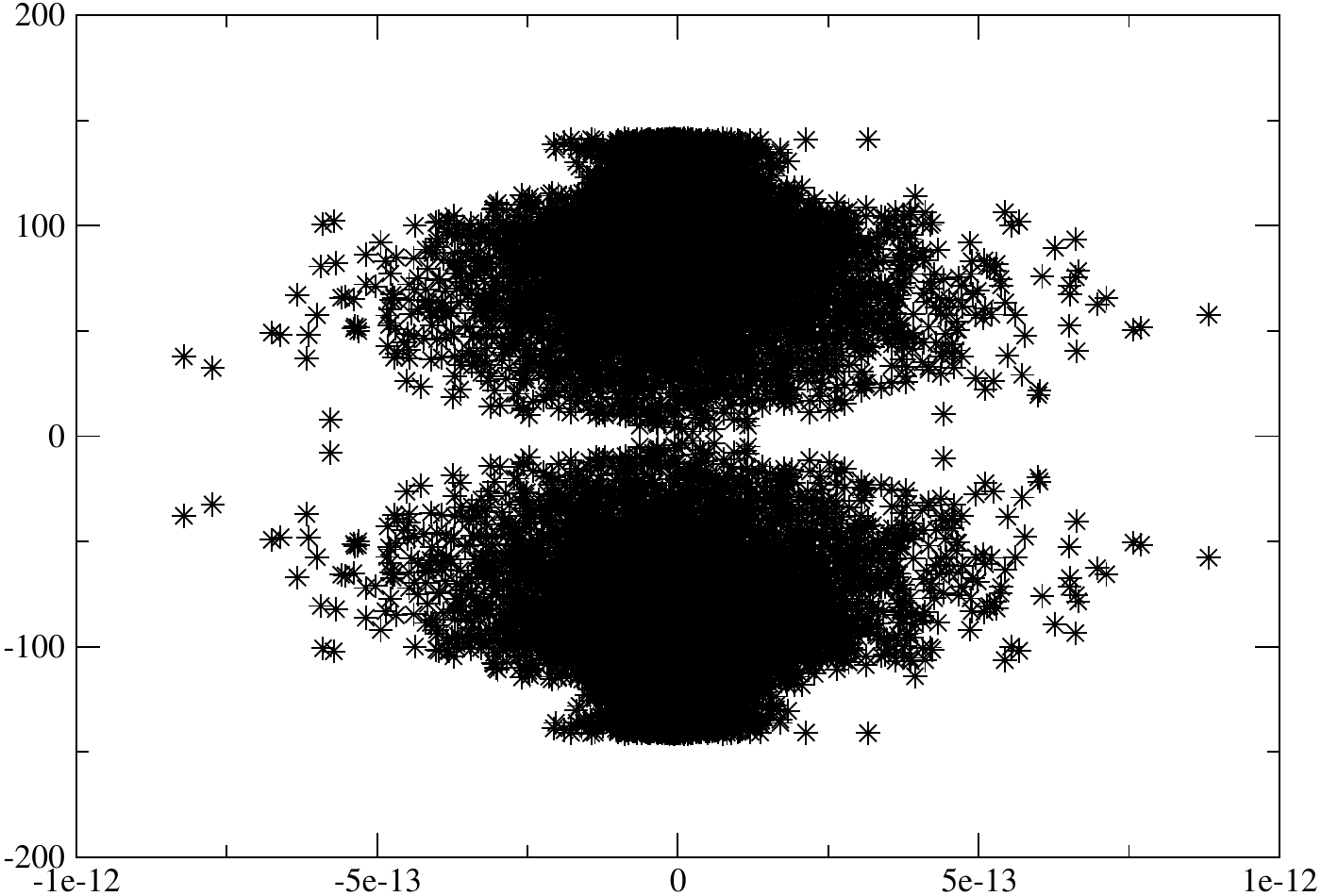}
  \caption{Homogeneous dimension 2 for $L=1$ and $N=100$}
  \label{dim2c0homogen}
\end{subfigure}
\hfill
 \begin{subfigure}{0.40\textwidth}
  \centering
  \includegraphics[width=\linewidth]{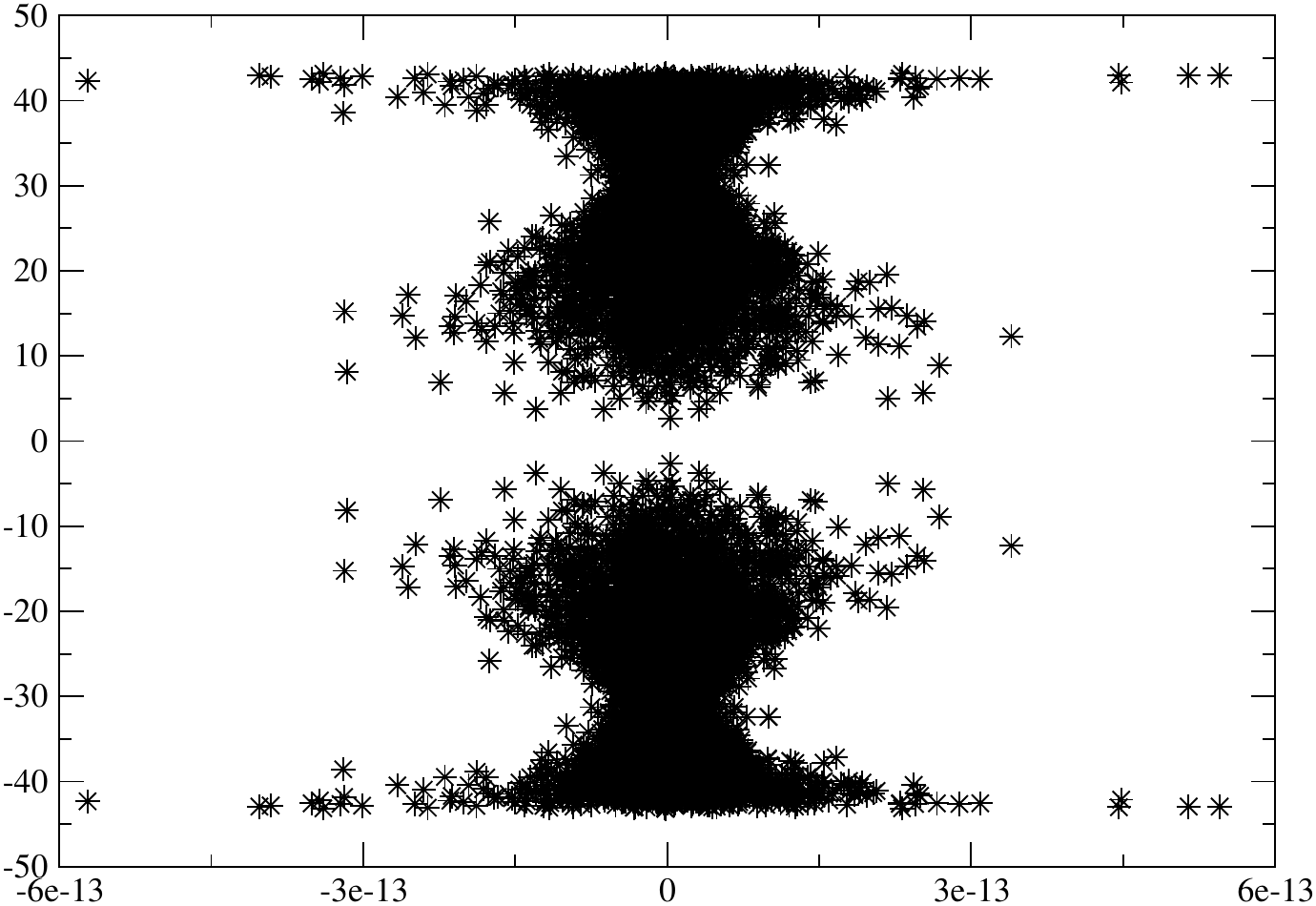}
  \caption{Homogeneous dimension 3 for $L=1$ and $N=25$}
  \label{dim3c0homogen}
\end{subfigure}
\caption{Homogeneous boundary conditions for $n=2$, $3$ and $c=0$}
\label{coefc0dim23}
\end{figure}

\begin{figure}[hbt!]
\centering
 \begin{subfigure}{0.40\textwidth}
  \centering
  \includegraphics[width=\linewidth]{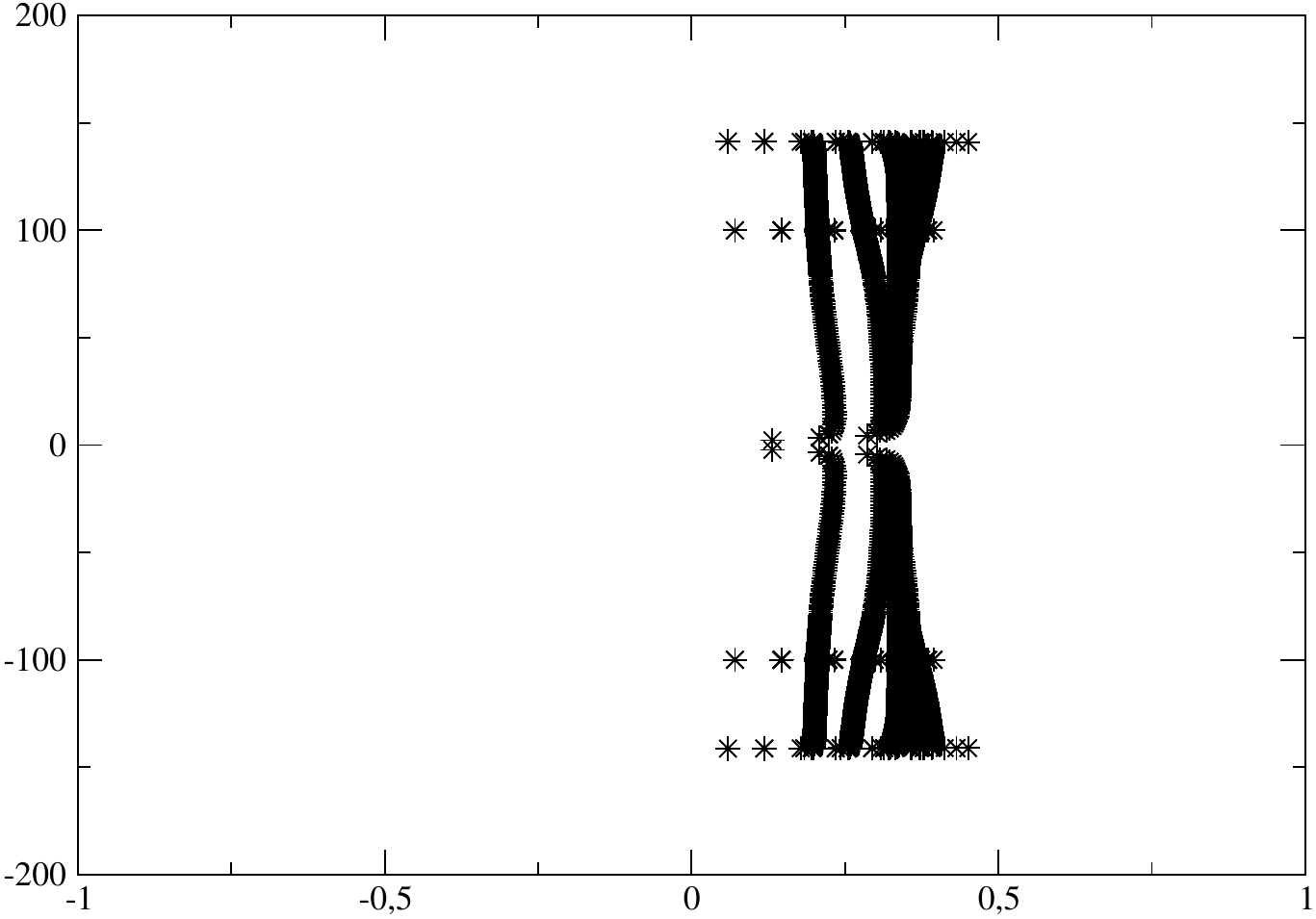}
  \caption{Homogeneous dimension 2 for $L=1$ and $N=100$}
  \label{dim2c1homogen}
\end{subfigure}
\hfill
 \begin{subfigure}{0.40\textwidth}
  \centering
  \includegraphics[width=\linewidth]{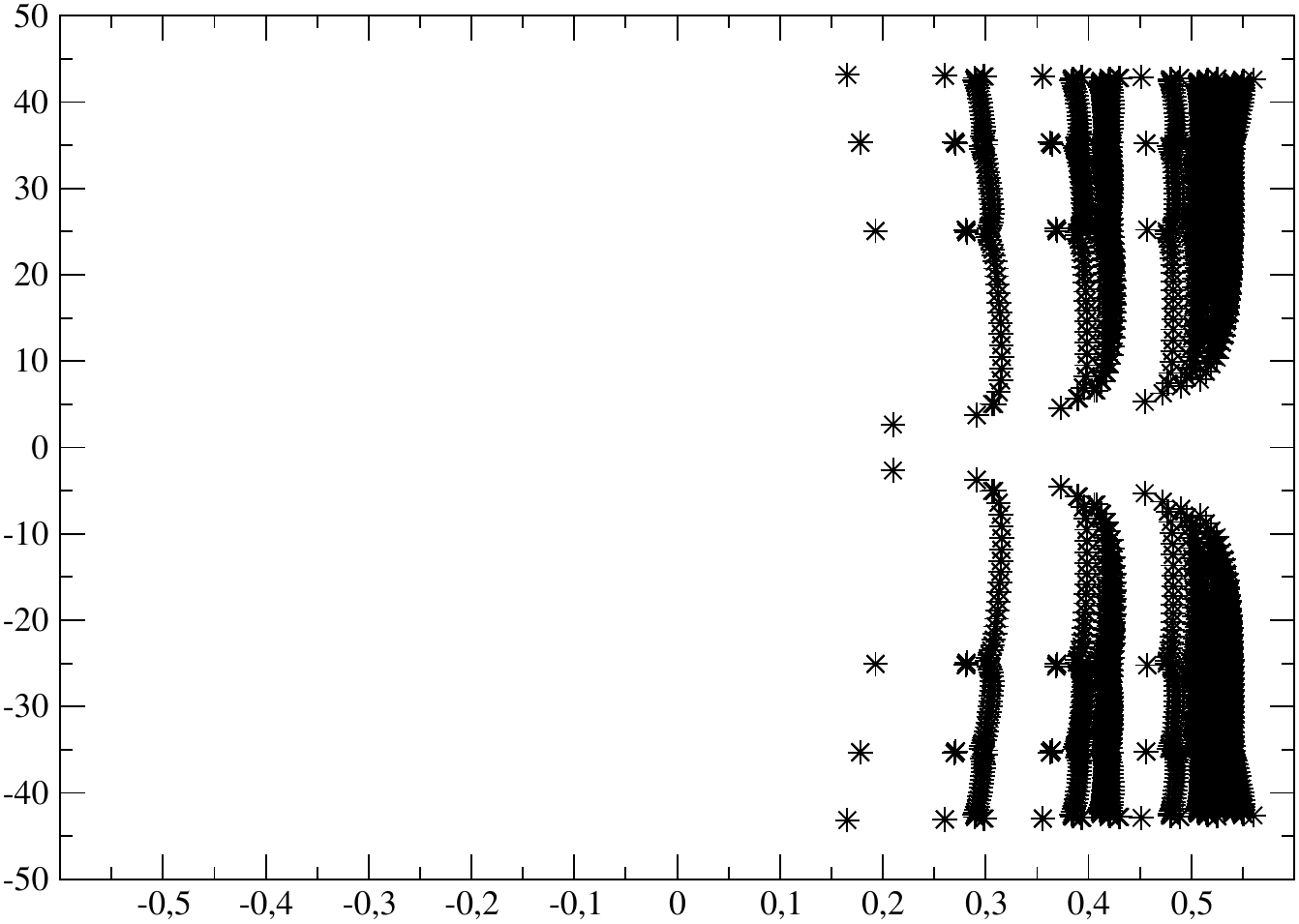}
  \caption{Homogeneous dimension 3 for $L=1$ and $N=25$}
  \label{dim3c1homogen}
\end{subfigure}
\caption{Homogeneous boundary conditions for $n=2$, $3$ and $c=0.5$}
\label{coefc1dim23}
\end{figure}

\begin{figure}[hbt!]
\centering
 \begin{subfigure}{0.40\textwidth}
  \centering
  \includegraphics[width=\linewidth]{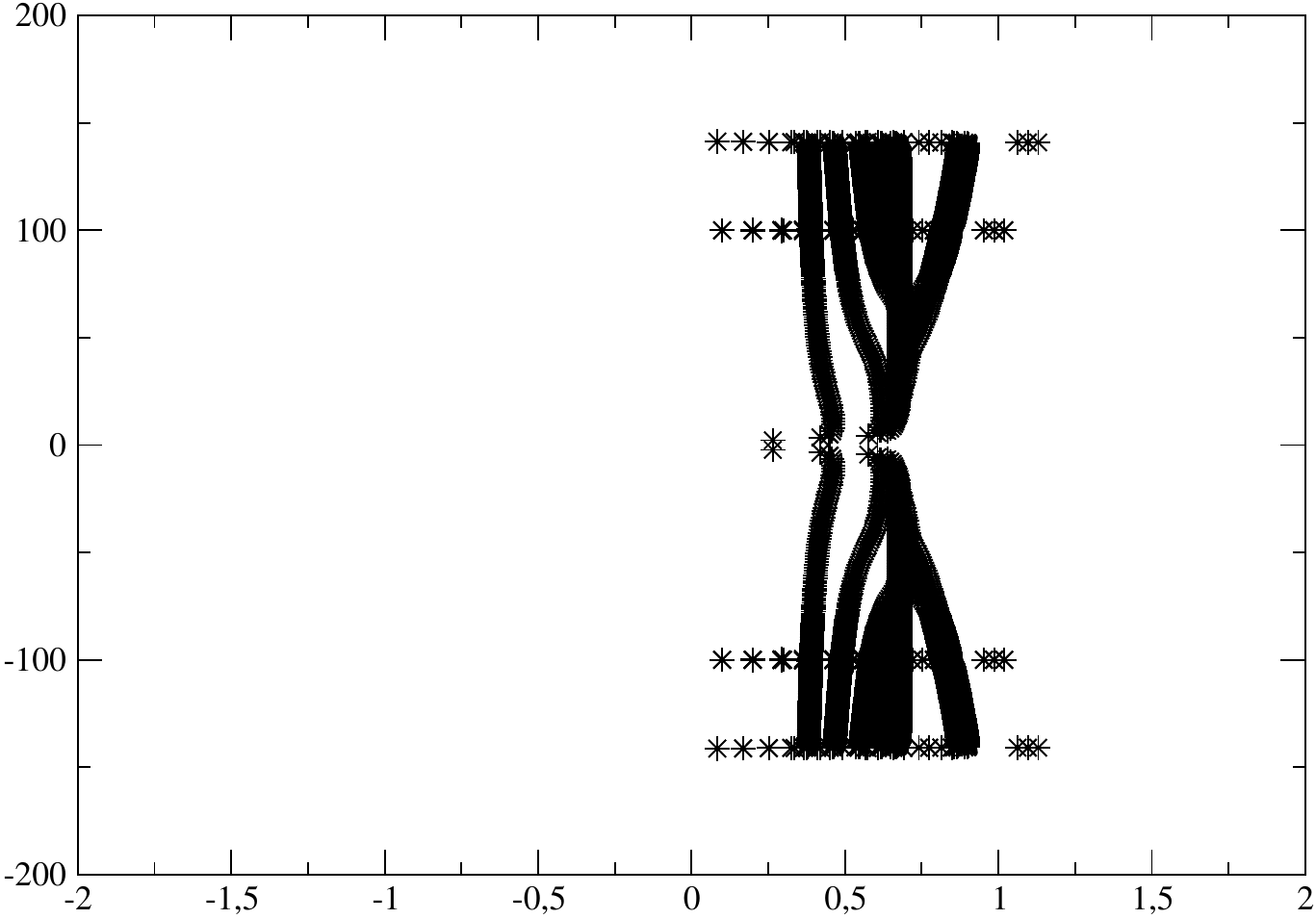}
  \caption{Homogeneous dimension 2 for $L=1$ and $N=100$}
  \label{dim2c2homogen}
\end{subfigure}
\hfill
 \begin{subfigure}{0.40\textwidth}
  \centering
  \includegraphics[width=\linewidth]{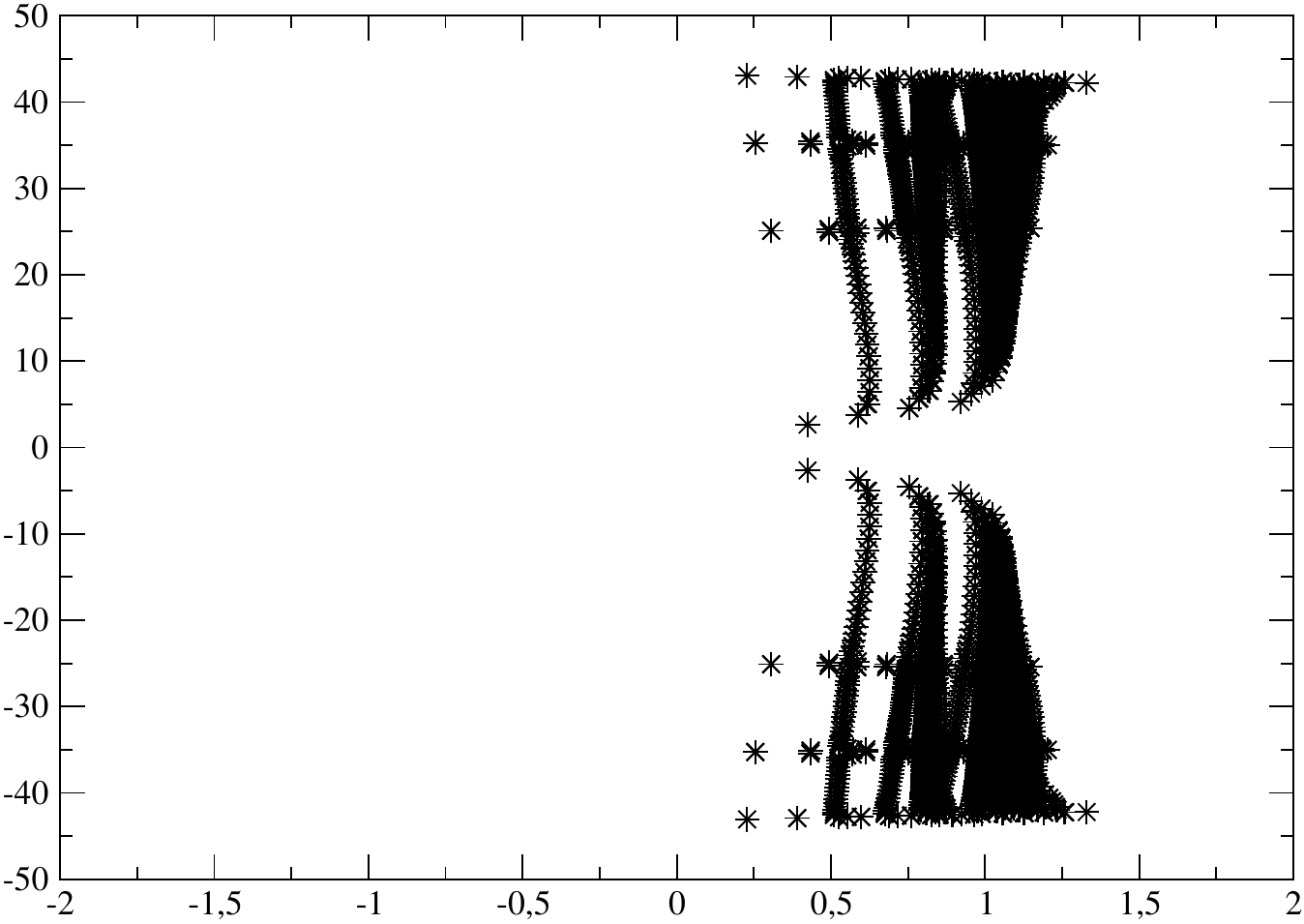}
  \caption{Homogeneous dimension 3 for $L=1$ and $N=25$}
  \label{dim3c2homogen}
\end{subfigure}
\caption{Homogeneous boundary conditions for $n=2$, $3$ and $c=1$}
\label{coefc2dim23}
\end{figure}

\begin{figure}[hbt!]
\centering
 \begin{subfigure}{0.40\textwidth}
  \centering
  \includegraphics[width=\linewidth]{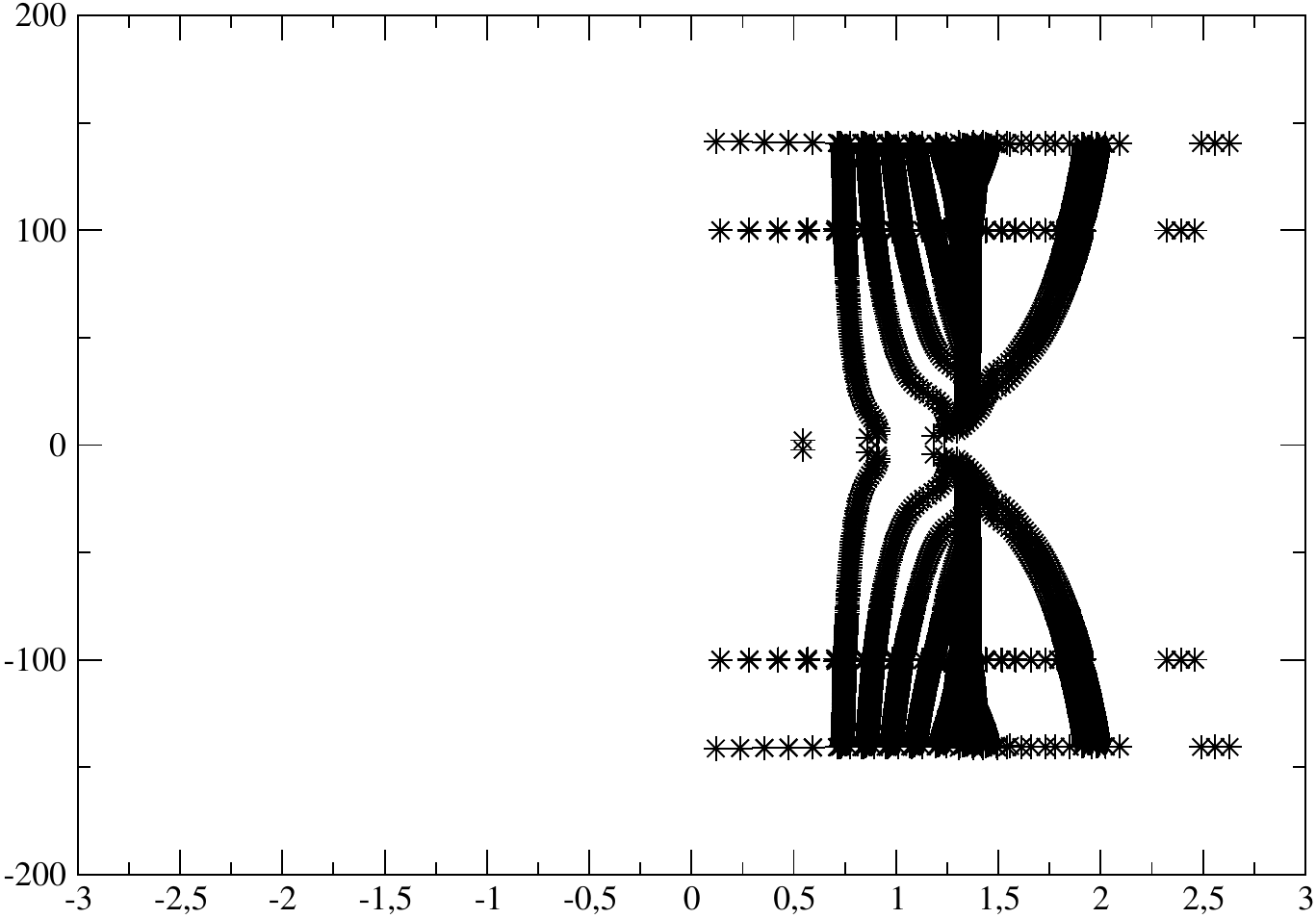}
  \caption{Homogeneous dimension 2 for $L=1$ and $N=100$}
  \label{dim2c4homogen}
\end{subfigure}
\hfill
 \begin{subfigure}{0.40\textwidth}
  \centering
  \includegraphics[width=\linewidth]{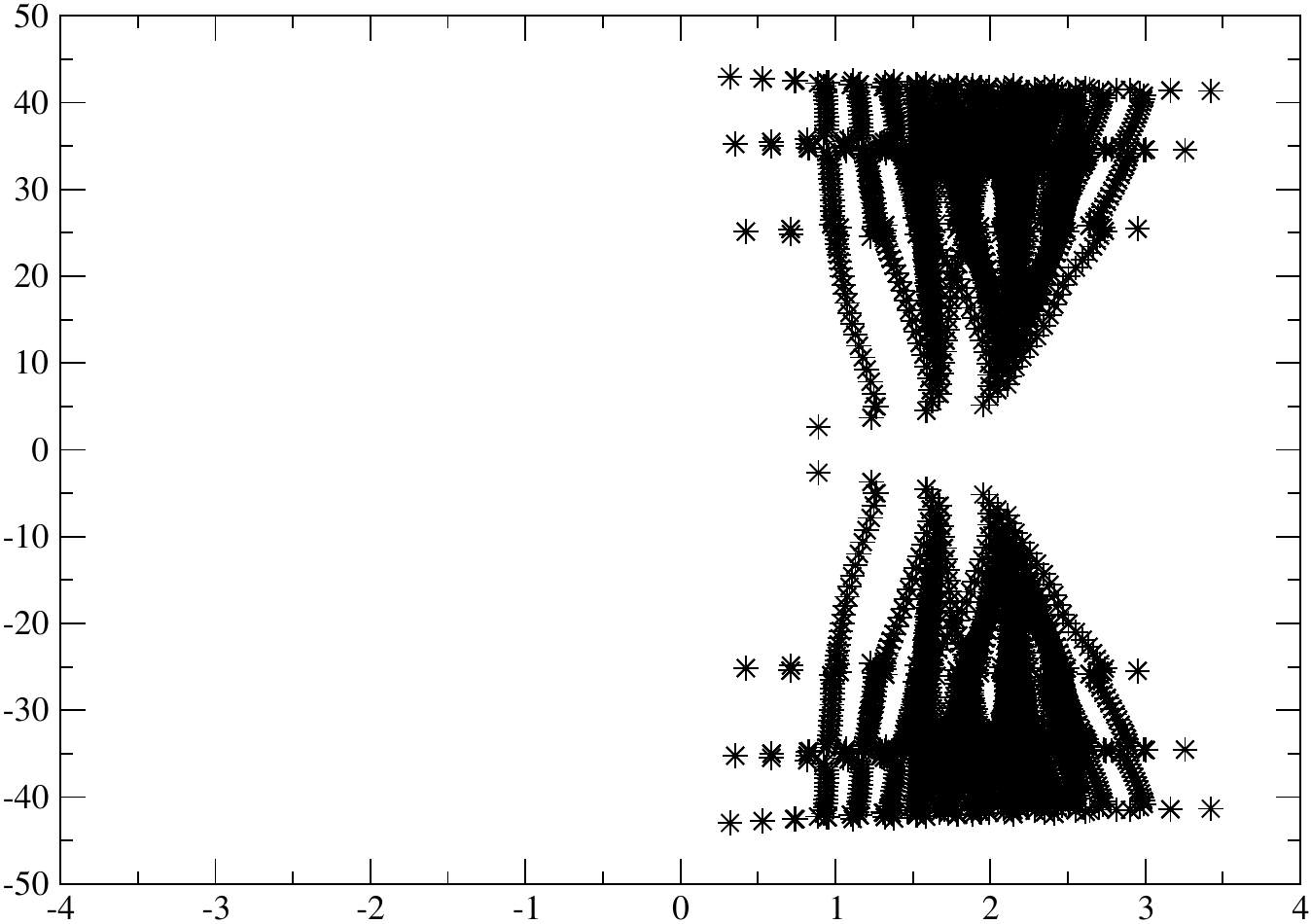}
  \caption{Homogeneous dimension 3 for $L=1$ and $N=25$}
  \label{dim3c4homogen}
\end{subfigure}
\caption{Homogeneous boundary conditions for $n=2$, $3$ and $c=2$}
\label{coefc4dim23}
\end{figure}

\begin{figure}[hbt!]
\centering
 \begin{subfigure}{0.40\textwidth}
  \centering
  \includegraphics[width=\linewidth]{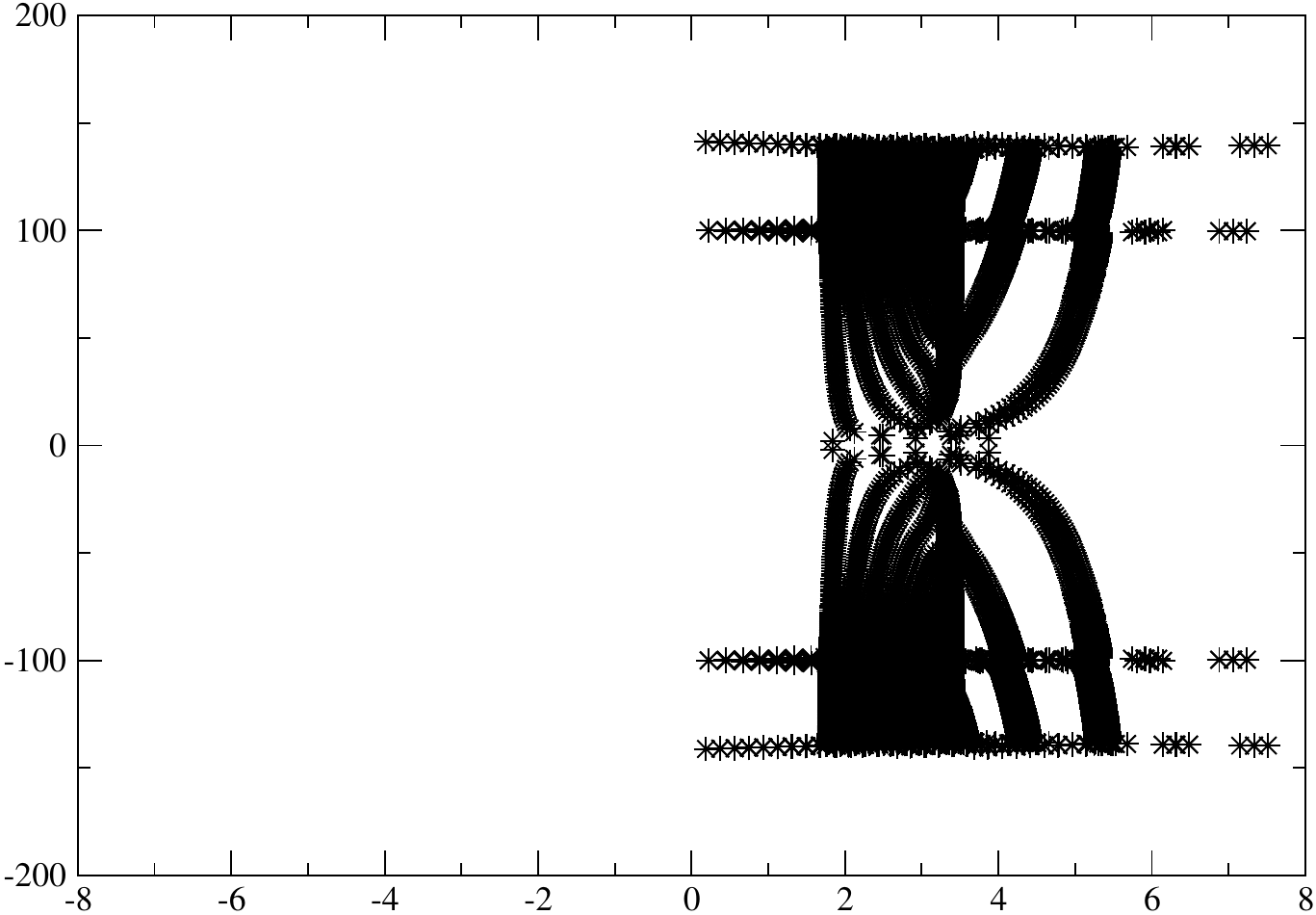}
  \caption{Homogeneous dimension 2 for $L=1$ and $N=100$}
  \label{dim2c10homogen}
\end{subfigure}
\hfill
 \begin{subfigure}{0.40\textwidth}
  \centering
  \includegraphics[width=\linewidth]{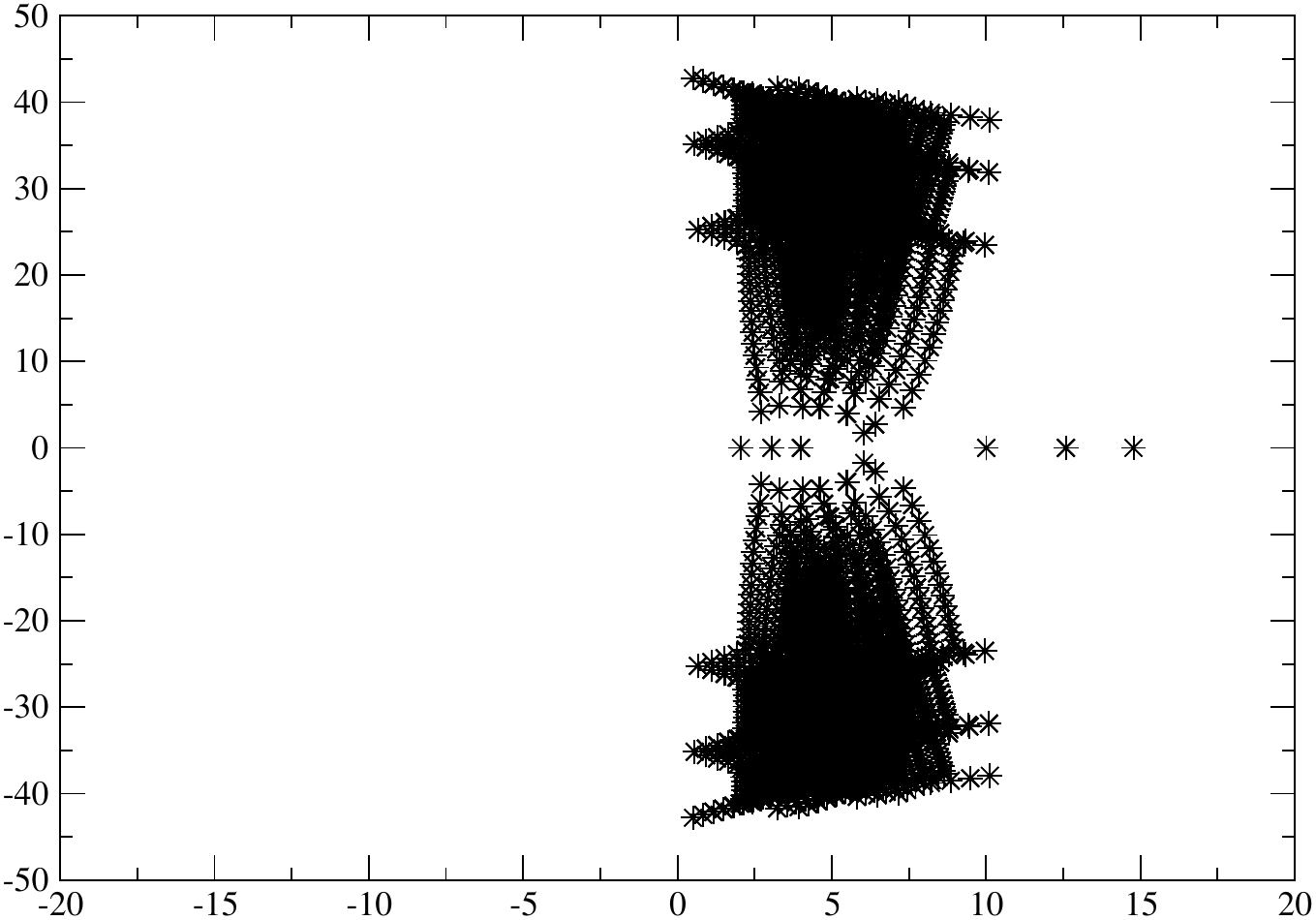}
  \caption{Homogeneous dimension 3 for $L=1$ and $N=25$}
  \label{dim3c10homogen}
\end{subfigure}
\caption{Homogeneous boundary conditions for $n=2$, $3$ and $c=5$}
\label{coefc10dim23}
\end{figure}

\begin{figure}[hbt!]
\centering
 \begin{subfigure}{0.40\textwidth}
  \centering
  \includegraphics[width=\linewidth]{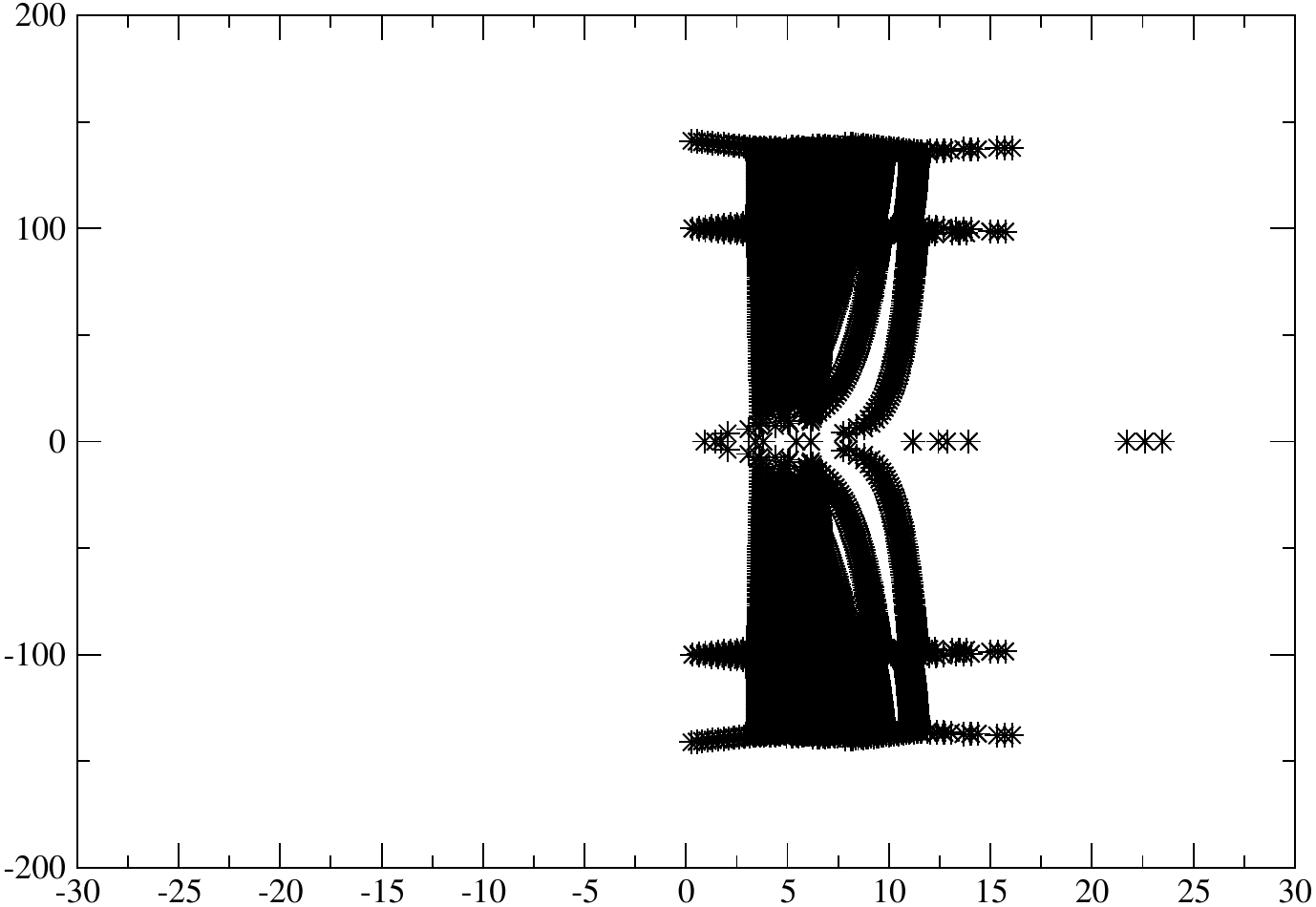}
  \caption{Homogeneous dimension 2 for $L=1$ and $N=100$}
  \label{dim2c20homogen}
\end{subfigure}
\hfill
 \begin{subfigure}{0.40\textwidth}
  \centering
  \includegraphics[width=\linewidth]{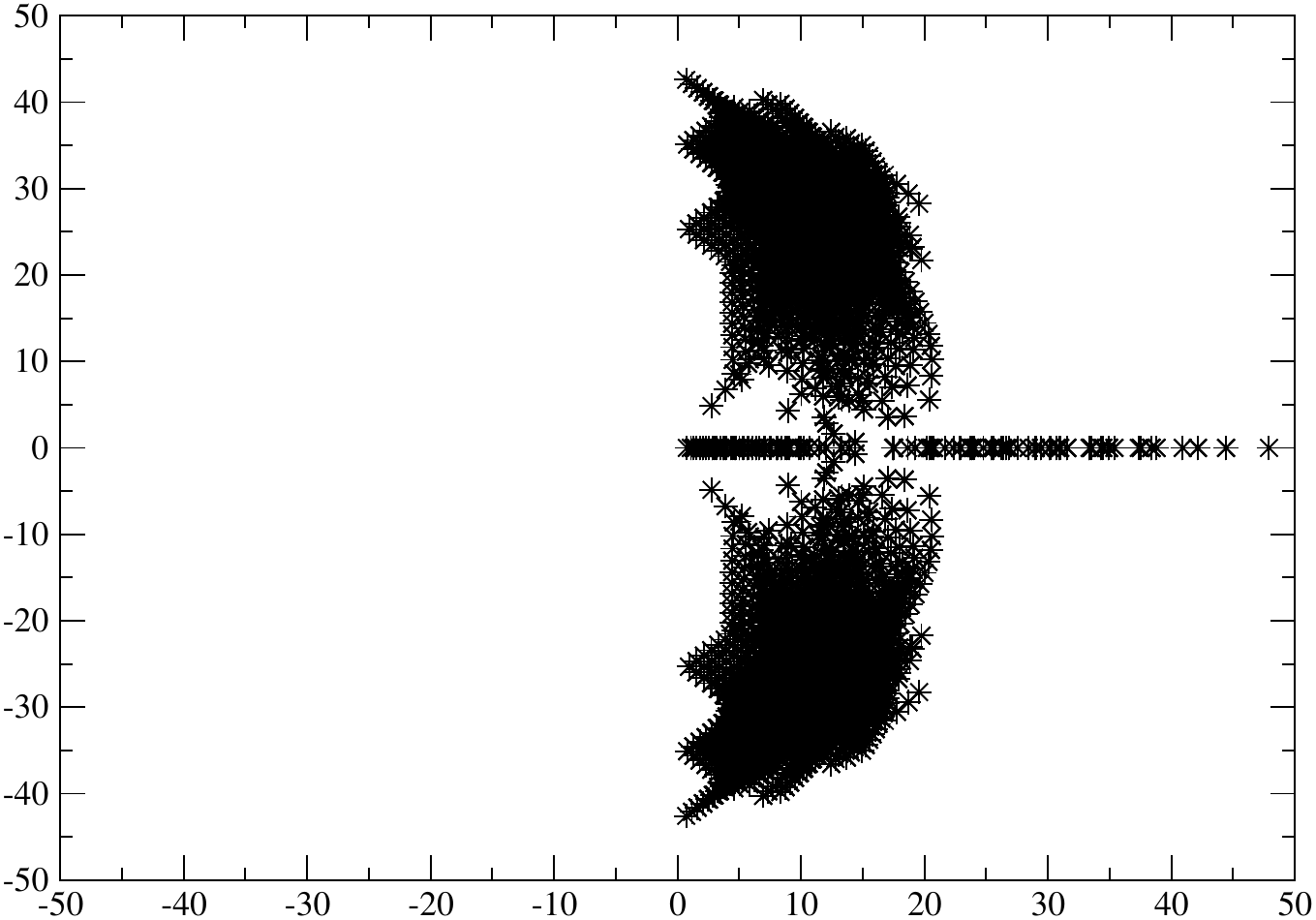}
  \caption{Homogeneous dimension 3 for $L=1$ and $N=25$}
  \label{dim3c20homogen}
\end{subfigure}
\caption{Homogeneous boundary conditions for $n=2$, $3$ and $c=10$}
\label{coefc20dim23}
\end{figure}

\begin{figure}[hbt!]
\centering
 \begin{subfigure}{0.40\textwidth}
  \centering
   \includegraphics[width=\linewidth]{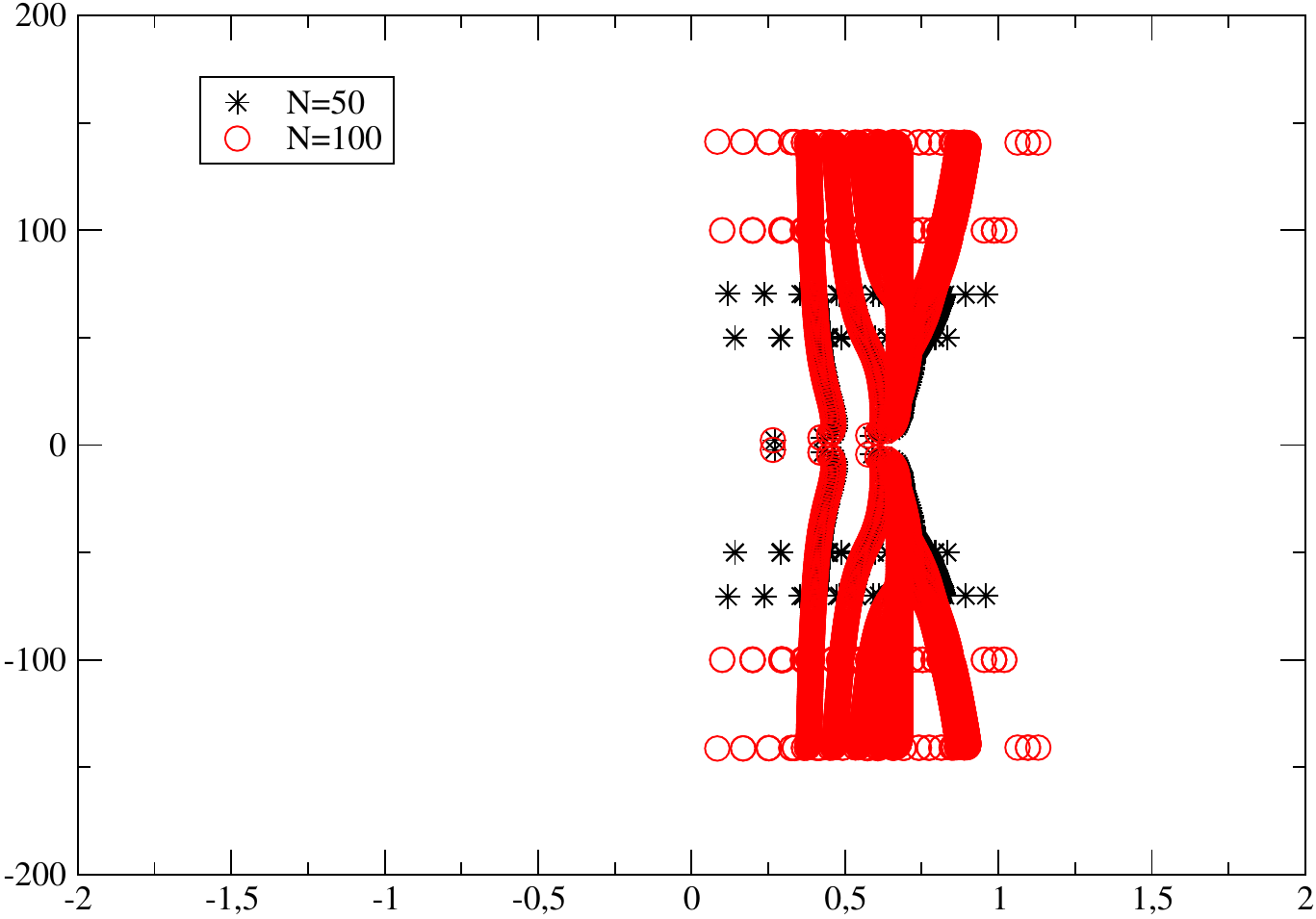}
  \caption{Comparison in dimension 2 with the same length of interval $L=1$ and a different number of partition points $N=50$ and $N=100$}
  \label{compare2Dsamediff}
\end{subfigure}
\hfill
 \begin{subfigure}{0.40\textwidth}
  \centering
  \includegraphics[width=\linewidth]{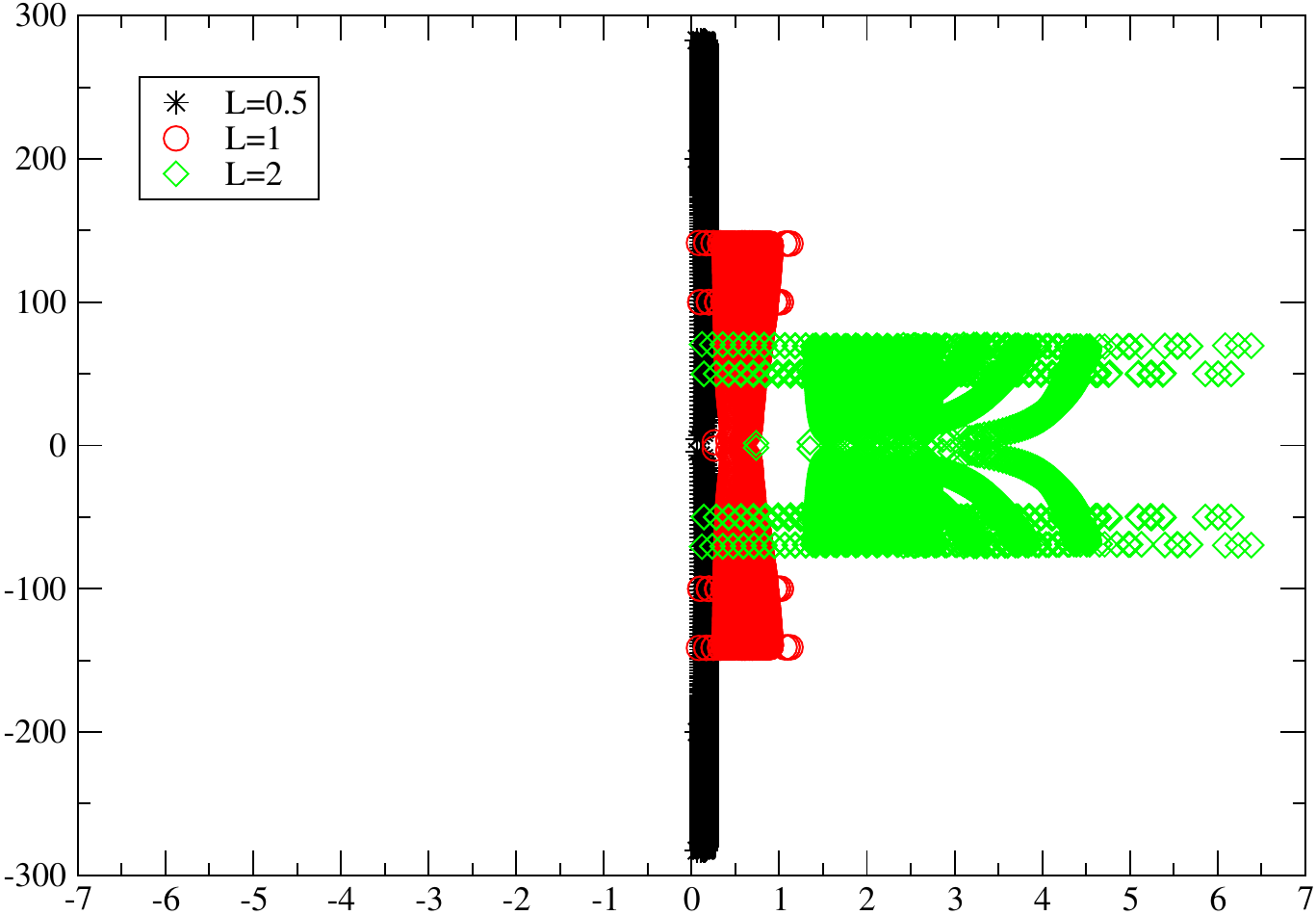}
  \caption{Comparison in dimension 2 with different lengths of interval $L=0.5$, $L=1$, $L=2$ and the same number of partition points $N=100$ }
  \label{compar2Ddiffsam}
\end{subfigure}
\caption{Homogenous boundary conditions for $n=2$ and $c=1$}
\label{comparc2dim3homogen}
\end{figure}

\begin{figure}[hbt!]
\centering
 \begin{subfigure}{0.40\textwidth}
  \centering
   \includegraphics[width=\linewidth]{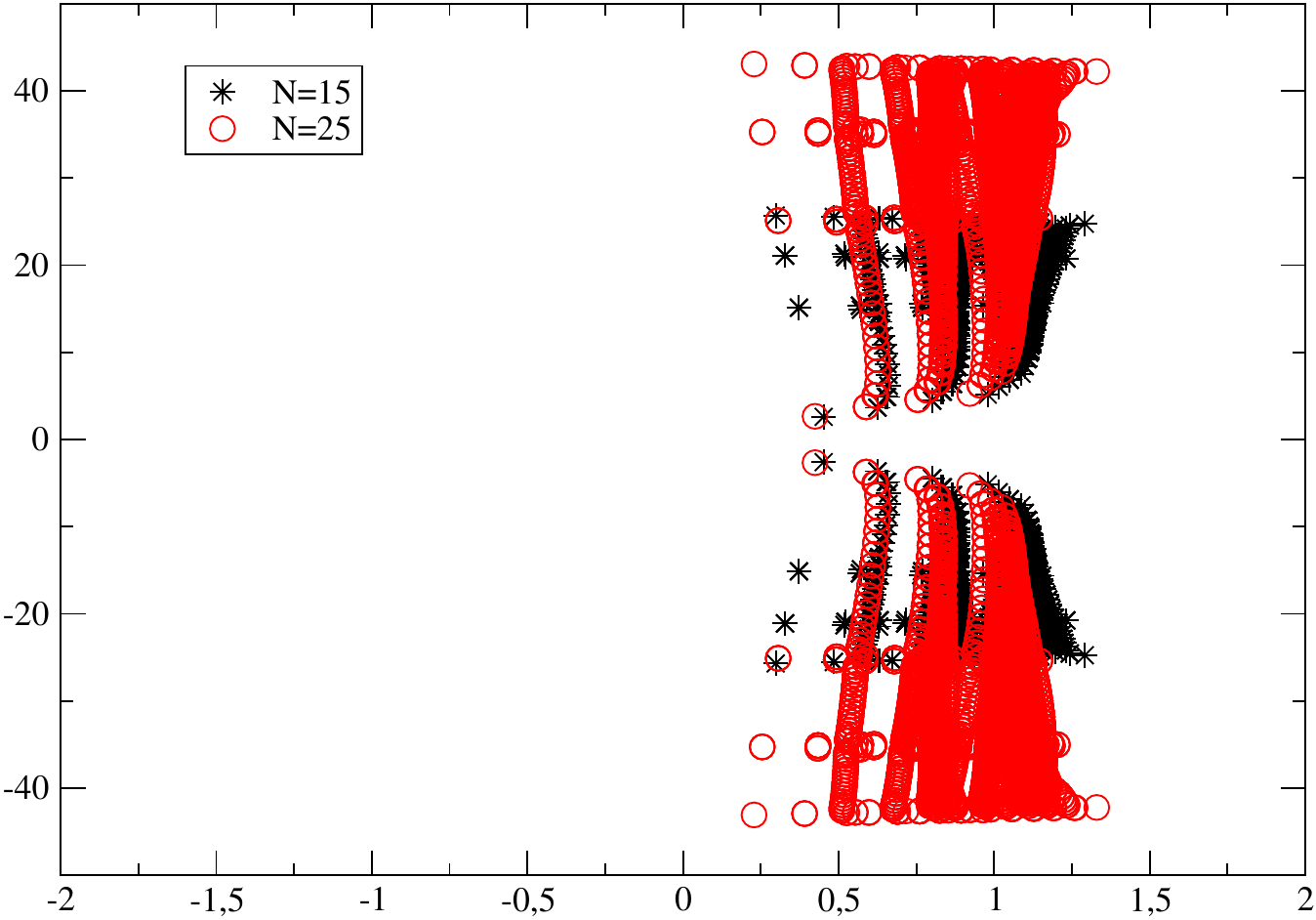}
  \caption{Comparison in dimension 3 with the same length of interval $L=1$ and a different number of partition points $N=15$ and $N=25$}
  \label{compare3Dsamediff}
\end{subfigure}
\hfill
 \begin{subfigure}{0.40\textwidth}
  \centering
  \includegraphics[width=\linewidth]{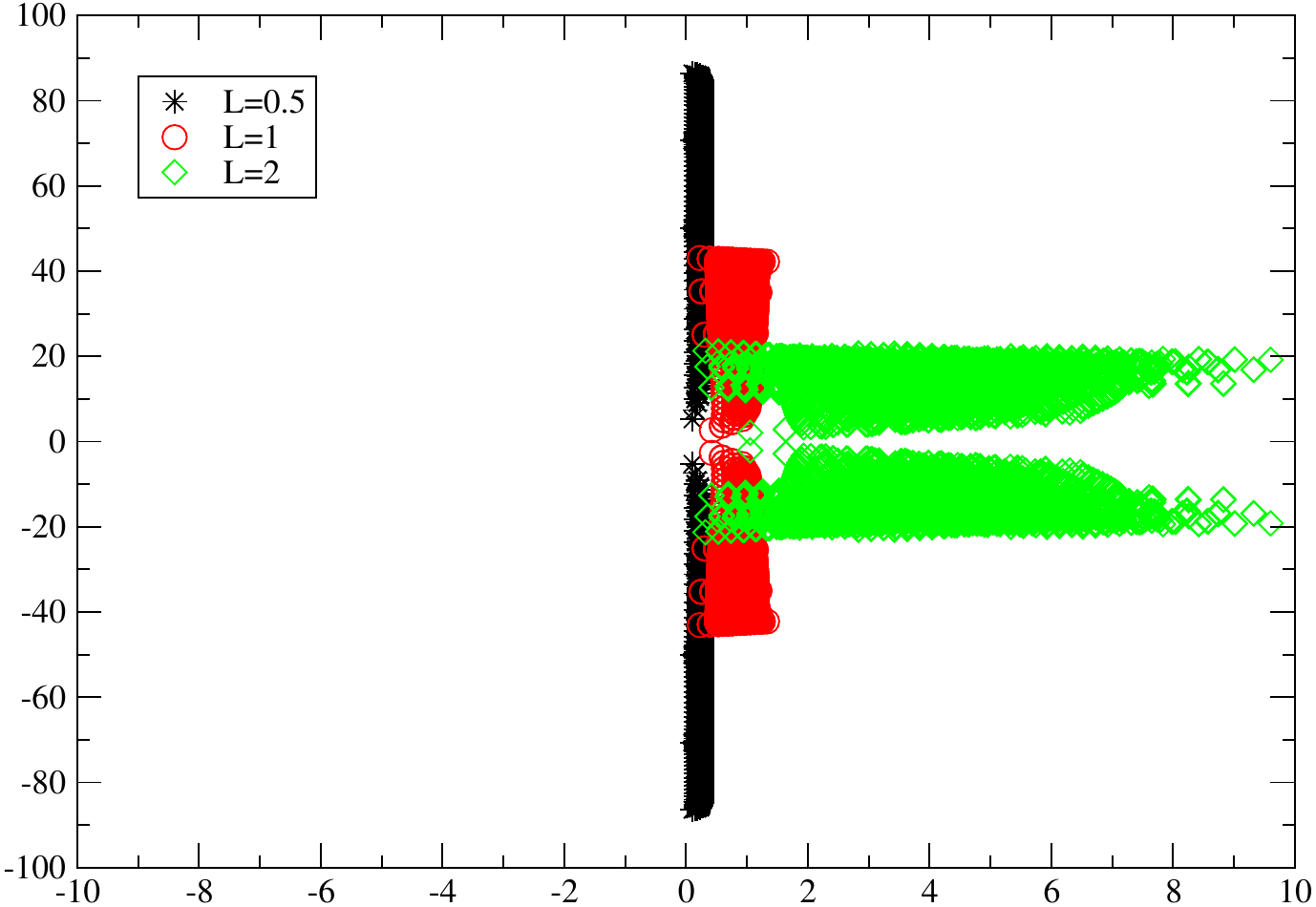}
  \caption{Comparisons in dimension 3 with different lengths of interval $L=0.5$, $L=1$, $L=2$ and the same number of partition points $N=25$ }
  \label{compar3Ddiffsam}
\end{subfigure}
\caption{Homogenous boundary conditions for $n=3$ and $c=1$}
\label{comparc2dim3homogen}
\end{figure}

\begin{figure}[hbt!]
\centering
 \begin{subfigure}{0.40\textwidth}
  \centering
\includegraphics[width=\linewidth]{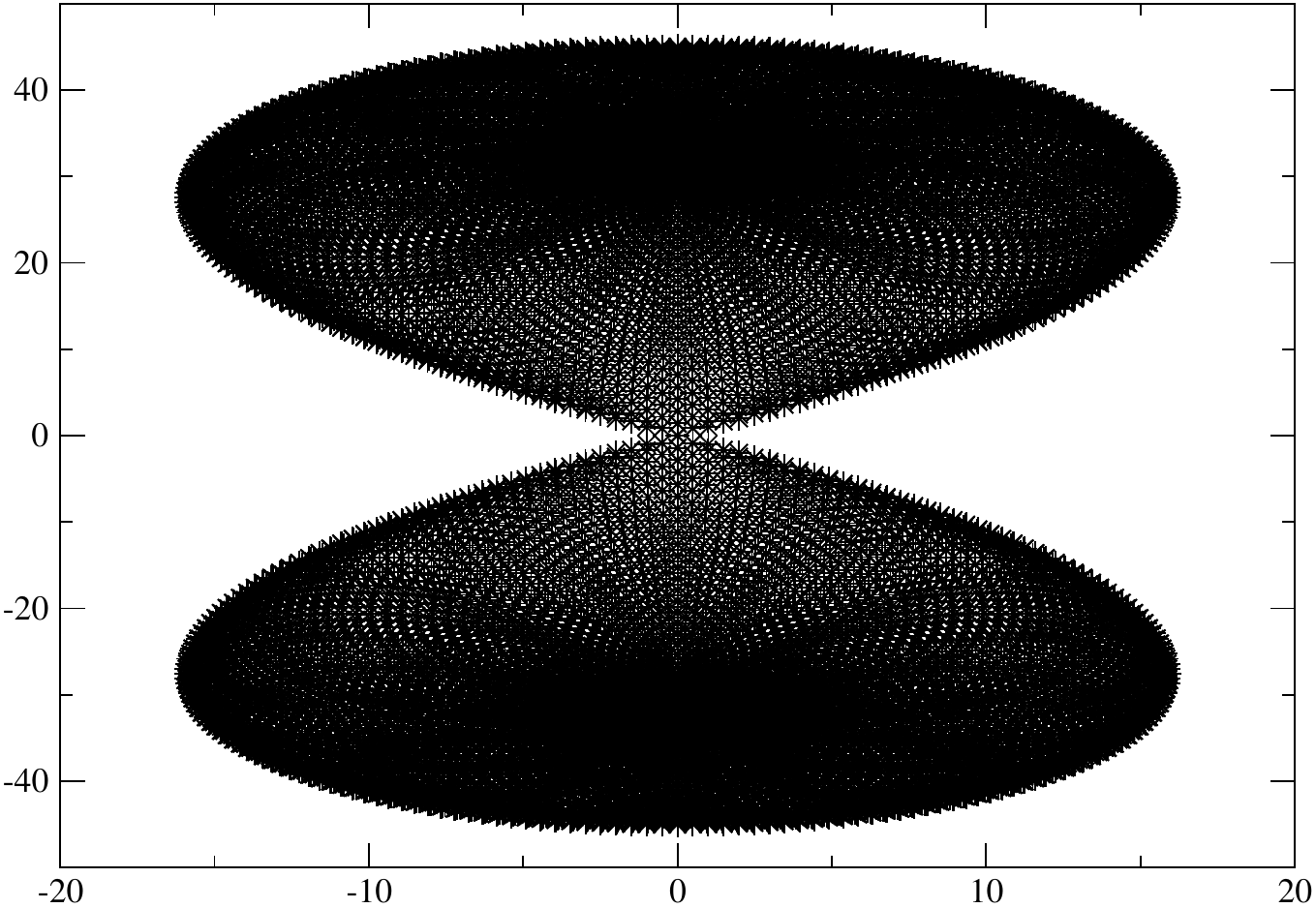}
\caption{$a_1(x_1)=a_2(x_2)=1$}
\label{dim2periodic11}
\end{subfigure}
\hfill
 \begin{subfigure}{0.40\textwidth}
  \centering
\includegraphics[width=\linewidth]{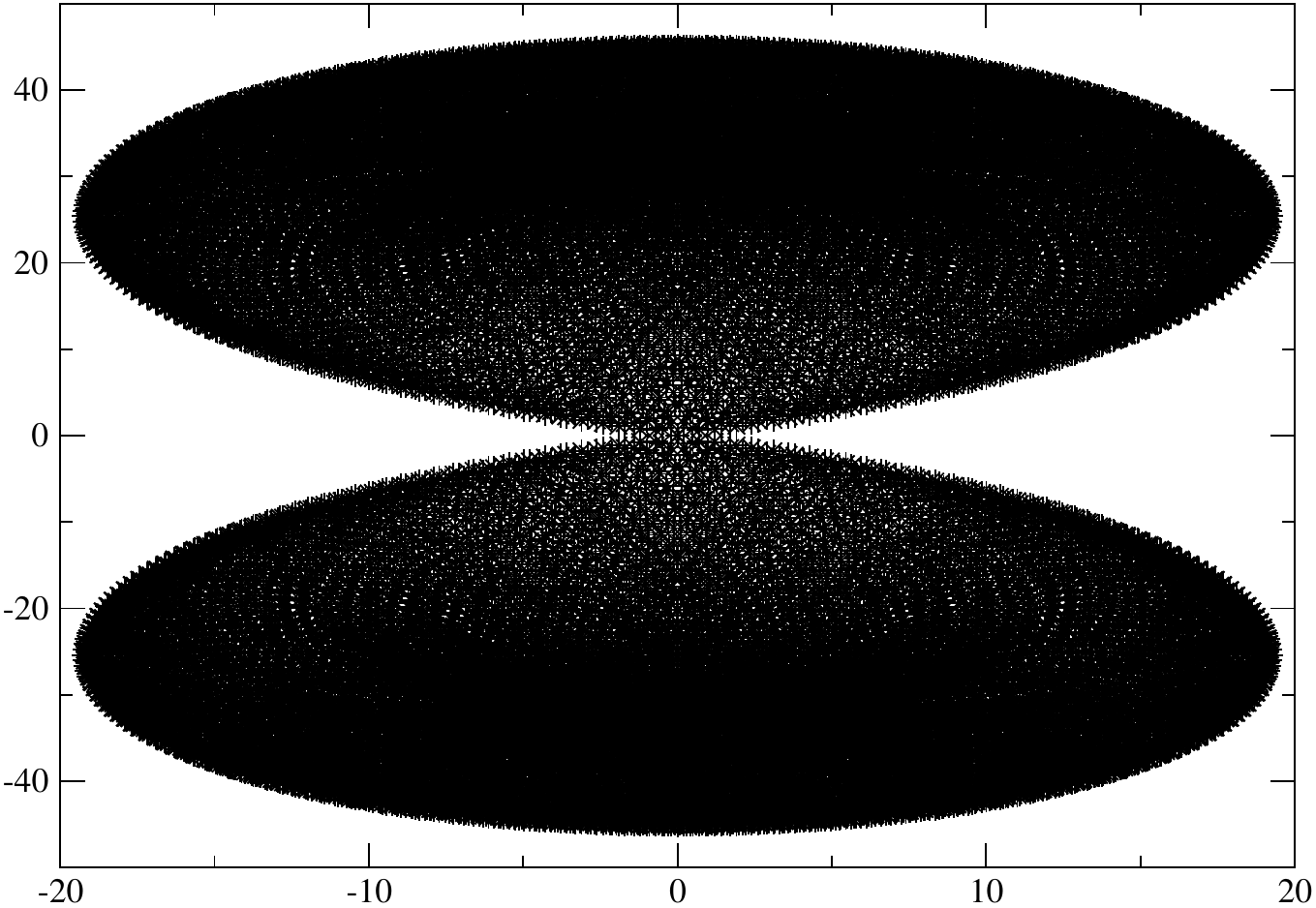}
\caption{$a_1(x_1)=1$ and $a_2(x_2)=\sqrt{2}$}
\label{dim2periodic1sqrt2}
\end{subfigure}
\caption{Periodic boundary conditions for $n=2$}
\end{figure}

\begin{figure}[hbt!]
\centering
 \begin{subfigure}{0.40\textwidth}
  \centering
\includegraphics[width=\linewidth]{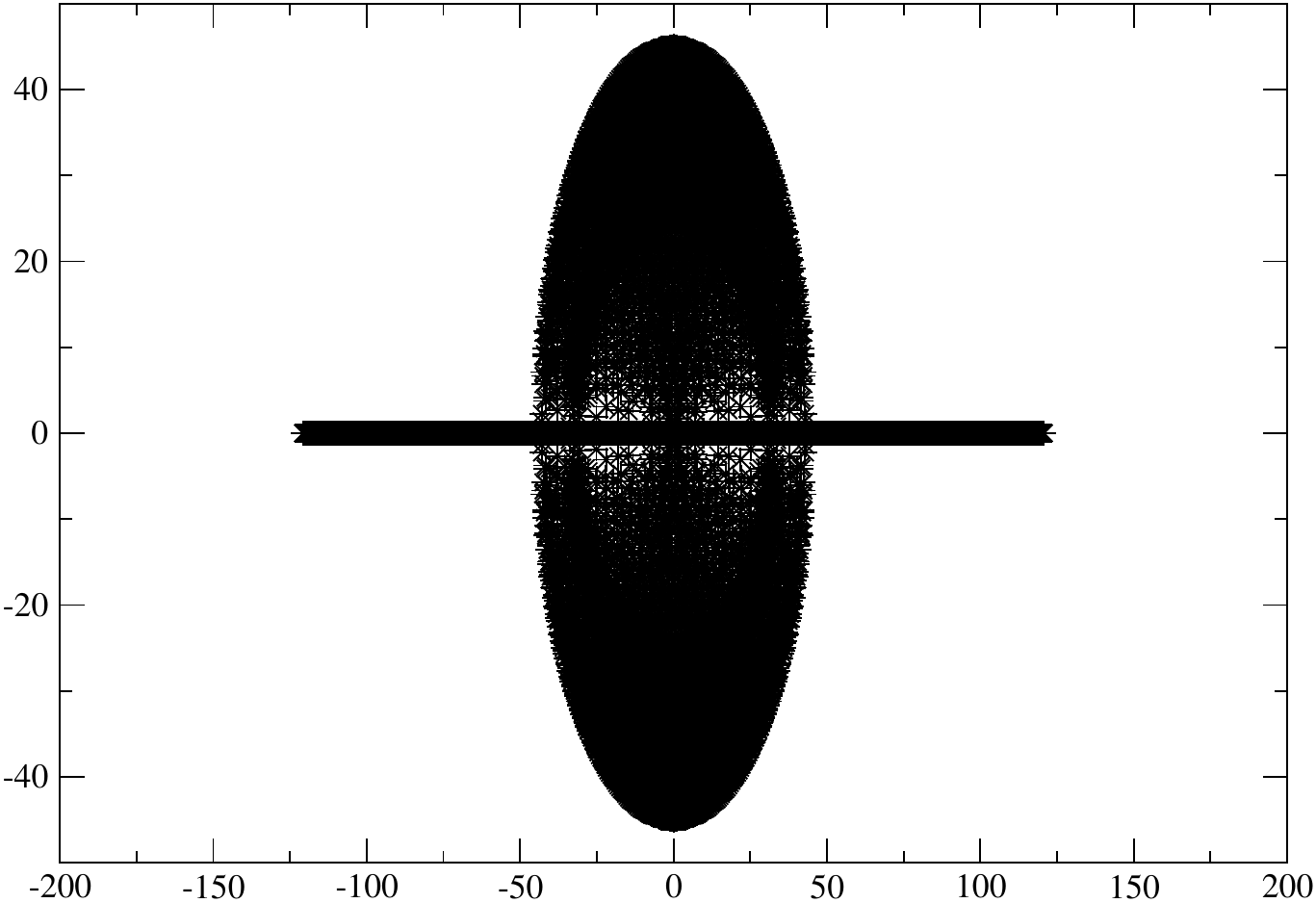}
\caption{$a_1(x_1)=1$ and $a_2(x_2)=5\sqrt{2}$}
\label{dim2periodic15sqrt2}
\end{subfigure}
\hfill
 \begin{subfigure}{0.40\textwidth}
  \centering
\includegraphics[width=\linewidth]{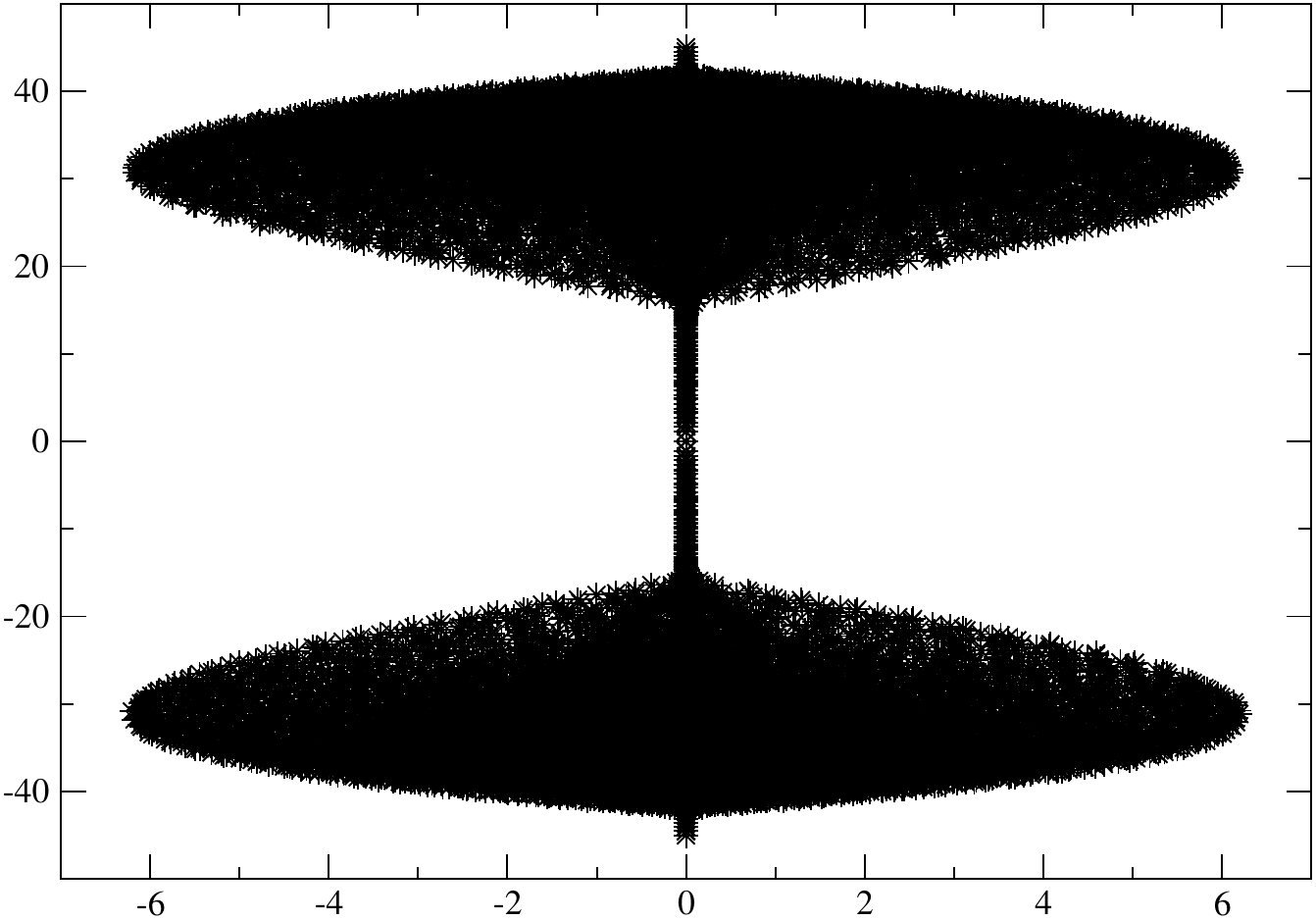}
\caption{$a_1(x_1)=\sin(x_1)$ and $a_2(x_2)=\sin(x_2)$}
\label{dim2periodicsinsin}
\end{subfigure}
\caption{Periodic boundary conditions for $n=2$}
\end{figure}

\begin{figure}[hbt!]
\centering
 \begin{subfigure}{0.40\textwidth}
  \centering
\includegraphics[width=\linewidth]{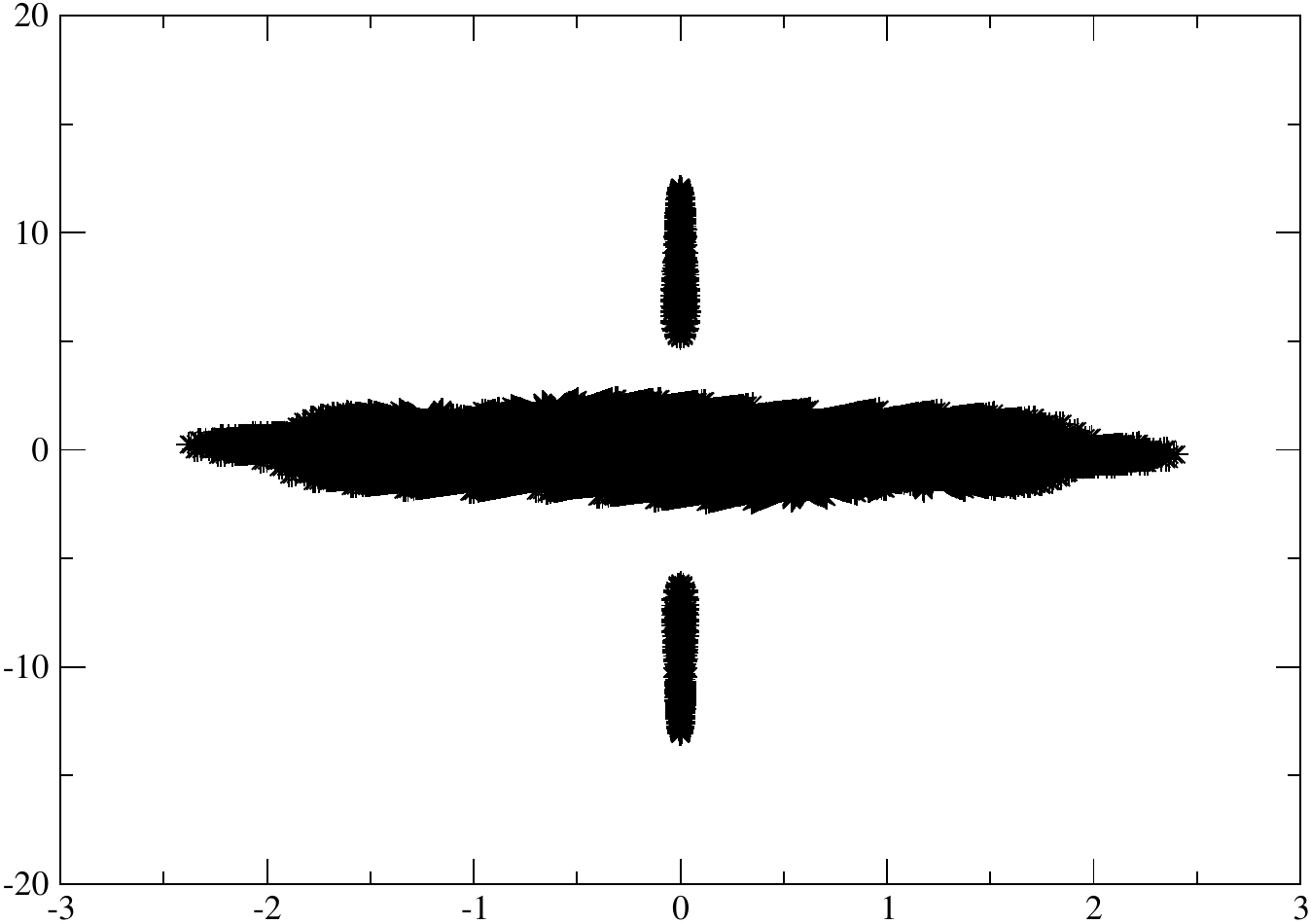}
\caption{$a_1(x_1)=a_2(x_2)=a_3(x_3)=1$}
\label{dim3periodic111}
\end{subfigure}
\hfill
 \begin{subfigure}{0.40\textwidth}
  \centering
\includegraphics[width=\linewidth]{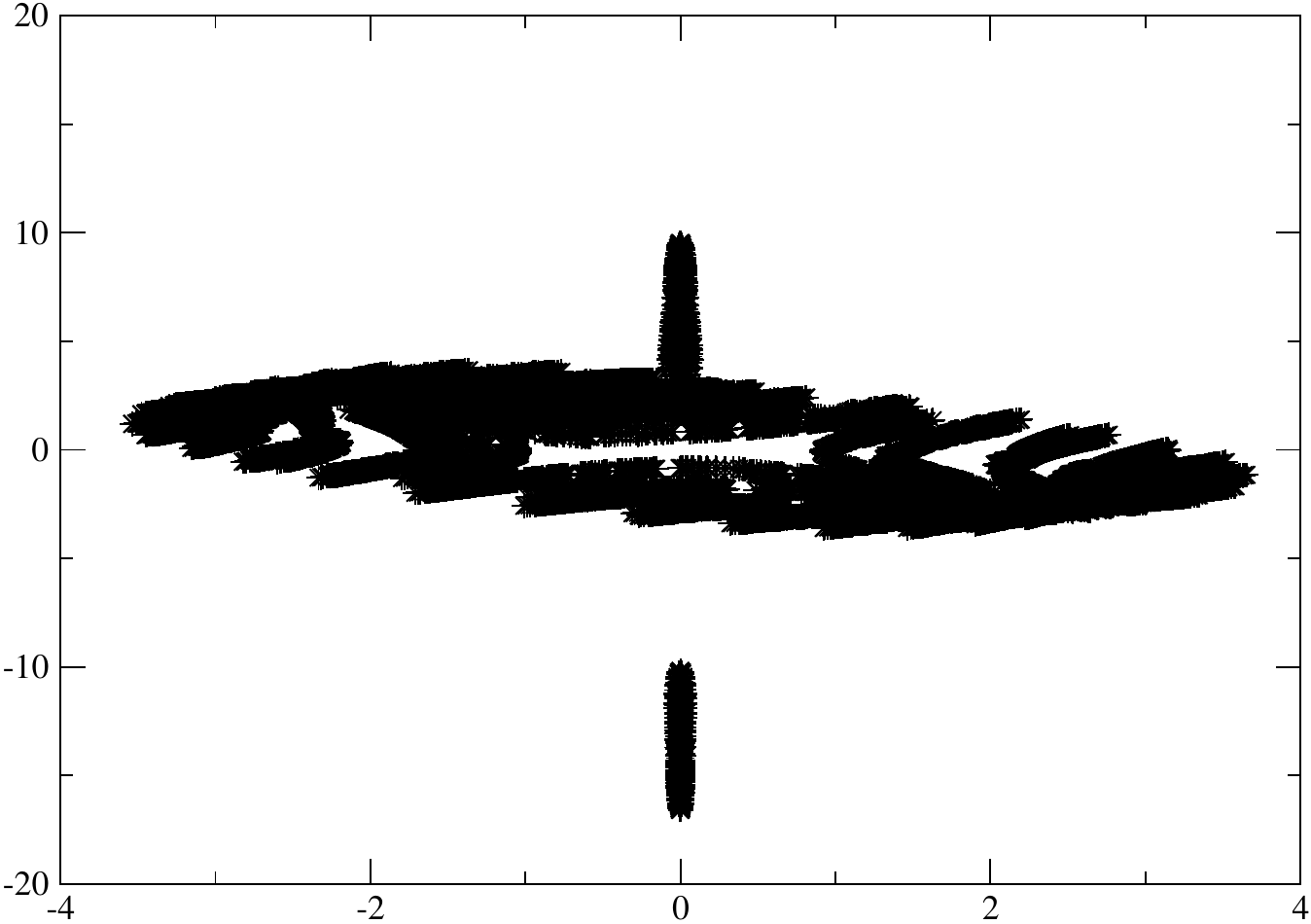}
\caption{$a_1(x_1)=1$, $a_2(x_2)=\sqrt{2}$ and $a_3(x_3)=5\sqrt{2}$}
\label{dim3periodic1sqrt25sqrt2}
\end{subfigure}
\caption{Periodic boundary conditions for $n=3$}
\end{figure}

\begin{figure}[hbt!]
\centering
 \begin{subfigure}{0.40\textwidth}
  \centering
\includegraphics[width=\linewidth]{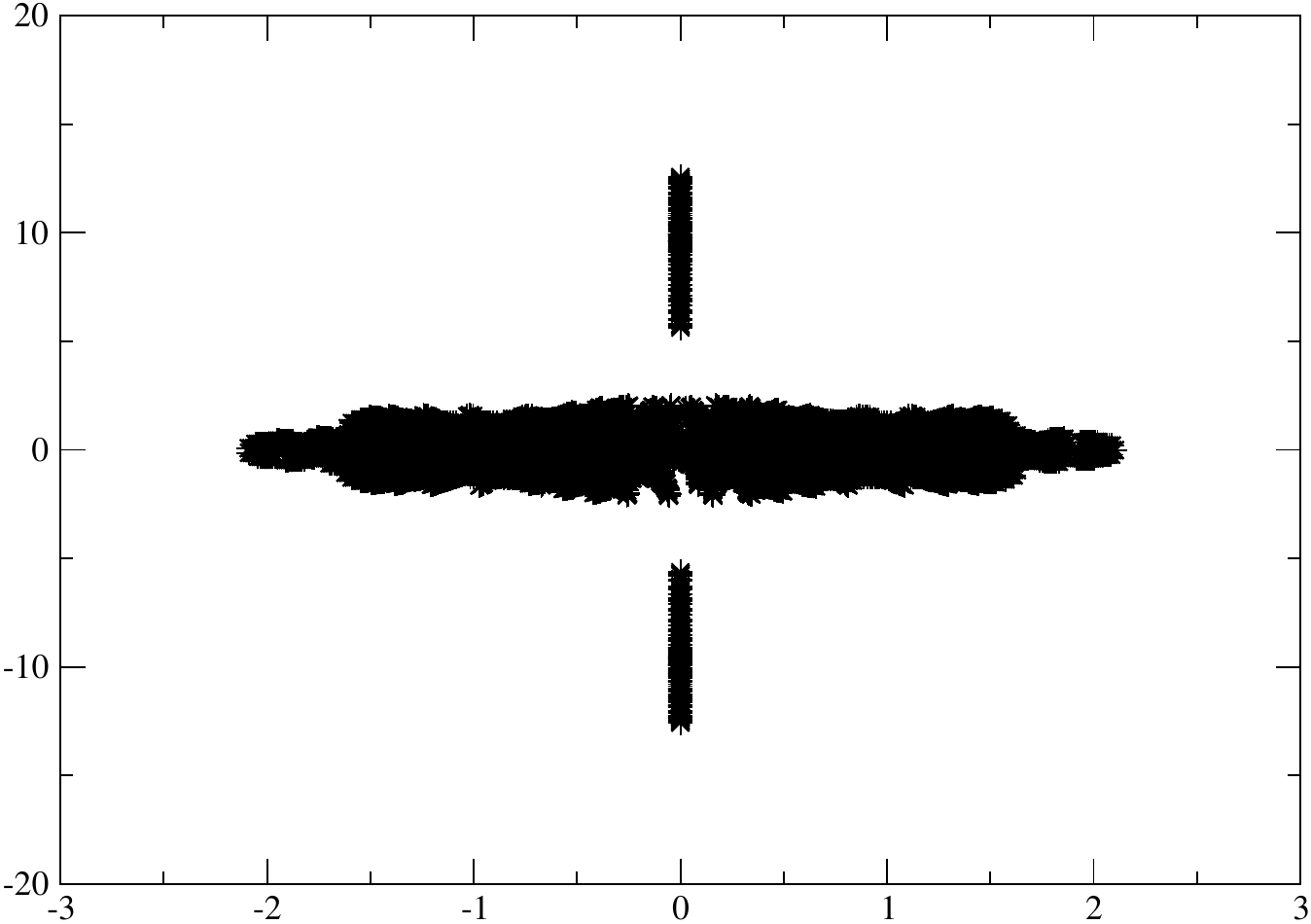}
\caption{$a_1(x_1)= \sin(x_1)$, $a_2(x_2)= \sin(x_2)$ and $a_3(x_3)=\sin(x_3)$}
\label{dim3periodicsinsinsin}
\end{subfigure}
\hfill
\caption{Periodic boundary conditions for $n=3$}
\end{figure}